\documentclass[a4paper,12pt,leqno]{amsart}
\usepackage{latexsym}
\usepackage[all]{xy}

\usepackage{amssymb} 
\usepackage{amsmath} 
\usepackage{amscd}

\def\P{{\mathbf{P}}}
\def\Z{{\mathbb{Z}}}
\def\K{{\mathbb{K}}}
\def\CC{{\mathbb{C}}}
\def\R{{\mathbb{R}}}
\def\C{{\mathbb{C}}}
\def\A{{\mathcal{A}}}
\def\B{{\mathcal{B}}}
\def\C{{\mathcal{C}}}
\def\DF{{\mathcal{DF}}}
\def\SF{{\mathcal{SF}}}
\def\PDC{{\mathcal{IPD}}}

\DeclareMathOperator{\codim}{codim}

\DeclareMathOperator{\coker}{coker}
\DeclareMathOperator{\Der}{Der}

\DeclareMathOperator{\pd}{pd}

\DeclareMathOperator{\POexp}{POexp}

\numberwithin{equation}{section}

\newcommand{\owari}{\hfill$\square$}

\theoremstyle{break}
\newtheorem{theorem}{Theorem}[section]
\newtheorem{prop}[theorem]{Proposition}
\newtheorem{cor}[theorem]{Corollary}
\newtheorem{lemma}[theorem]{Lemma}
\newtheorem{define}[theorem]{Definition}

\newtheorem{example}[theorem]{Example}
\newtheorem{problem}[theorem]{Problem}
\newtheorem{conj}[theorem]{Conjecture}

\title{Projective dimensions of hyperplane arrangements}
\author{Takuro Abe\address{Institute of Mathematics for Industry, Kyushu University, 
Motooka 744, Nishi-ku, 8190395 Fukuoka, Japan}\email{abe@imi.kyushu-u.ac.jp}
}
\date{\today}

\begin{document}

\maketitle

\begin{abstract}
We establish a general theory for projective dimensions of the logarithmic 
derivation modules of hyperplane arrangements. 
That includes the addition-deletion and restriction theorem, 
Yoshinaga-type result, and the division theorem for projective dimensions of hyperplane arrangements. 
They are generalizations of the
free arrangement cases, that can be regarded as the special case of our result when the projective 
dimension is zero.
The keys to prove them are several new methods to determine the surjectivity of the Euler and the Ziegler 
restriction maps, that is combinatorial when the projective dimension is not 
maximal for all localizations. Also, we introduce a new class of arrangements in which 
the projective dimension is comibinatorially determined. 
\end{abstract}

\section{Introduction}

\subsection{Setup and background}

Let $\K$ be an arbitrary field, $V=\K^\ell$, 
$S=\mbox{Sym}^*(V^*) \simeq \K[x_1,\ldots,x_\ell]$ and let $\Der S:=\oplus_{i=1}^\ell S\partial_{x_i}$ be the $S$-graded module of $\K$-linear $S$ derivations. 
Let $\A$ be an (central) 
\textbf{arrangement of hyperplanes} in $V$, i.e., a finite set of linear hyperplanes in $V$. 
For $H \in \A$, fix a linear form $\alpha_H$ such that $\ker \alpha_H=H$ and let 
$Q(\A):=\prod_{H \in \A} \alpha_H$. Then the \textbf{logarithmic derivation module $D(\A)$ of $\A$} is defined by 
$$
D(\A):=\{\theta \in \Der S \mid \theta(\alpha_H) \in S\alpha_H\ (\forall H \in \A)\}.
$$
$D(\A)$ is an $S$-graded reflexive module of rank $\ell$. 
We say that $\A$ is \textbf{free} with \textbf{exponents} $\exp(\A)=(d_1,\ldots,d_\ell)$ if 
there is a homogeneous $S$-basis $\theta_1,\ldots,\theta_\ell$ for $D(\A)$ such that $\deg \theta_i=d_i\ 
(i=1,\ldots,\ell)$. Here for $\theta \in \Der S$ we say that $\theta$ is \textbf{homogeneous of degree $d$} if 
$S \ni \theta(\alpha)=d$ for all $\alpha \in V^*$ with $\theta(\alpha) \neq 0$.

In the research of arrangements, the logarithmic derivation module has been one of the central topics. In particular, the freeness of $\A$ has been intensively studied from the beginning by K. Saito in \cite{Sa} showing that all Coxeter arrangements are free. The research of logarithmic derivation modules are from algebraic aspects, but recently several viewpoints are introduced to investigate them. The most important one is Terao's factorization theorem (Theorem \ref{Teraofactorization}). Let us explain for details.

Let $L(\A)$ be the intersection lattice (see Definition \ref{IL}) which remembers how hyperplanes in $\A$ intersect, i.e., it is the combinatorial structure of $\A$. From $L(\A)$ we can define the most important 
invariant $\chi(\A;t)=\sum_{i=0}^\ell (-1)^i b_i(\A)t^{\ell-i} $, the characteristic polynomial (Definition \ref{char}). By Brieskorn-Orlik-Solomon's result, we know that $b_i(\A)$ is the $i$-th Betti number of $V \setminus \cup_{H\in \A} H$ when $\K=\CC$. Thus $\chi(\A;t)$ is both combinatorial and topological 
invariants. Terao's factorization theorem asserts that if $\A$ is free with $\exp(\A)=(d_1,\ldots,d_\ell)$, then 
$$
\chi(\A;t)=\prod_{i=1}^\ell (t-d_i).
$$
Hence freeness is related to combinatorics, topology and algebraic geometry. Also recently, it was shown in \cite{AHMMS} that freeness plays the key rote to show the K\"{a}hler package and Poincar\`{e} duality of the regular nilpotent Hessenberg variety. Thus it is expected that the role played by the freeness will be more important in several research areas. 

In the study of free arrangements $\A$ and the logarithmic derivation modules $D(\A)$, 
the most useful result is the following.

\begin{theorem}[Terao's addition-deletion theorem, \cite{T1}]
Let $H \in \A$, $\A':=\A\setminus \{H\}$ and $\A^H:=\{H \cap L \mid 
L \in \A'\}$. Then two of the following three imply the 
third:
\begin{itemize}
\item[(1)]
$\A$ is free with $\exp(\A)=(d_1,\ldots,d_{\ell-1},d_\ell)$.
\item[(2)]
$\A'$ is free with $\exp(\A')=(d_1,\ldots,d_{\ell-1},d_\ell-1)$.
\item[(3)]
$\A^H$ is free with $\exp(\A^H)=(d_1,\ldots,d_{\ell-1})$.
\end{itemize}
In particular, all the three hold true if both $\A$ and $\A'$ are free.
\label{adddel}
\end{theorem}

Theorem \ref{adddel} is used to check the (non-)freeness of $\A$, and to construct free arrangements. 
Though Theorem \ref{adddel} was shown 40 years before, it is most used in these days in the freeness 
research. 
One point in Theorem \ref{adddel} is that to apply it we need two algebraic information, i.e., 
freeness and exponents. Let us re-consider Theorem \ref{adddel} by revising these two information follwoing the argument in \cite{A2}. First, by using Terao's factorization (Theorem \ref{Teraofactorization}), 
exponents can be replaced by 
combinatorial one, e.g., 
$$
\chi(\A^H;t) \mid \chi(\A;t) \iff \exp(\A^H) \subset \exp(\A)
$$
if $\A$ and $\A^H$ are free. Thus we can formulate Theorem \ref{adddel} in terms of freeness and 
the division of characteristic polynomials. In fact, to apply Theorem \ref{adddel}, 
the \textbf{$b_2$-equality} is sufficient (see Theorem \ref{division} for details on the $b_2$-equality):
$$
b_2(\A)=b_2(\A^H)+|\A^H|(|\A|-|\A^H|).
$$
We say that the \textbf{$b_2$-equality holds for $(\A,H)$} when 
the above holds. 
Then it is easy to show that 
\begin{eqnarray*}
&\ &\exp(\A^H) \subset \exp(\A)\ \mbox{or}\ 
\exp(\A^H) \subset \exp(\A')\ \mbox{or}\  
|\exp(\A') \cap \exp(\A) |=\ell-1\\
&\Rightarrow& \chi(\A^H;t) \mid \chi(\A;t) \\
&\Rightarrow& 
b_2(\A)=b_2(\A^H)+|\A^H|(|\A|-|\A^H|).
\end{eqnarray*}
Here we used the famous \textbf{deletion-restriction formula}:
$$
\chi(\A;t)=\chi(\A';t)-\chi(\A^H;t).
$$
Thus we can replace inclusions of exponents 
by the division of characterisric polynomials or 
the $b_2$-equality. 
The advantage to use the $b_2$-equality instead of the exponents is that we can assume it for non-free arrangements on which we cannot define exponents.
In fact, by using Theorems \ref{adddel} and \ref{division}, we can obtain the following formulation.

\begin{cor}[Addition-deletion theorem based on the $b_2$-equality]
Assume that the $b_2$-equality holds for $(\A,H)$. Then 
two of the following three imply the 
third:
\begin{itemize}
\item[(1)]
$\pd \A=0$.
\item[(2)]
$\pd \A'=0$.
\item[(3)]
$\pd \A^H=0$.
\end{itemize}
\label{adddel2}
\end{cor}

Here 
$
\pd \A:=\pd_S D(\A).
$
Clearly $\A$ is free if and only if $\pd \A=0$. Note that $\pd \A^H=\pd_{S^H} D(\A^H)$ for the coordinate ring 
$S^H$ of $H$.
Thus Theorem \ref{adddel} can be regarded as the result to \textbf{compare projective dimensions of 
logarithmic modules $D(\A),\ D(\A')$ and $D(\A^H)$ when at least two of them are zero, and 
the $b_2$-equality holds}. 

Also, it is easy to show that $\pd \A_X=0$ if $\pd \A=0$ (see Lemma \ref{pdlow} 
for example) for the localization $\A_X:=\{H \in \A \mid X \subset H\}$. Thus the assumption includes 
the information on the projective dimension on $\A_X$. We revise these two, i.e., \begin{itemize}
\item[(1)] 
The $b_2$-equality, or related local information in codimension three, and 
\item[(2)]
$\pd \A_X$ for $X \in L(\A^H)$.
\end{itemize}
to approach \textbf{the most important conjecture in the research of free arrangements} 
as follows:

\begin{conj}[Terao's conjecture]
Whether $\pd \A=0$ or not
depends only on $L(\A)$. Namely, if $\pd \A=0$, and there is $\B$ such that 
$L(\A)$ and $L(\B)$ 
are isomorphic as lattices, then $\pd \B=0$.
\label{TC}
\end{conj}

Since $ 0 \le \pd \A \le \ell-2$, Terao's conjecture 
is true if $\ell \le 2$. 
However, if $\ell \ge 3$, then almost 
nothing is known about Terao's conjecture. 
In fact, only few have been known to be (non-)combinatorial. For example, the ring structure of the cohomology ring of 
the complement of $\A$ is combinatorial when $\K=\CC$, but the fundamental group is not.

Now from our new viewpoint of freeness by $\pd \A$ and the 
$b_2$-equality, there is a new approach to Terao's conjecture. Namely, \textbf{establish the theory to 
control 
$\pd \A$ in terms of combinatoris}.  
In the rest of the subsections we introduce 
main results in this article. 


\subsection{NMPD and addition-deletion theorem for projetive dimensions}
Recall that $D(\A)$ is reflexive. Thus by the Auslander-Buchsbaum formula, it holds that 
$$
0 \le \pd \A:=\pd_S D(\A)
\le \ell-2.
$$
So let us give names when $\pd \A$ is the largest.  

\begin{define}
(1)\,\, 
We say that 
$\pd \A$ is \textbf{maximal} (or \textbf{$\A$ is of maximal projective dimension}) 
if $\pd  \A=\ell-2$. 

\noindent
(2)\,\,
We say that 
$\A$ is \textbf{not of maximal projective dimension (NMPD) along $H$} if 
$\pd \A_X < \codim  X-2$ 
for all $X \in L_{\ge 3}(\A^H)$, i.e., 
$\pd_{S_X} \A_X^e$ is not maximal for all $X \in L_{\ge 3}(\A^H)$. Here $S_X$ is the coordinate ring of 
$V/X$, and $\A_X^e$ is the essentialization of $\A_X$ (see 
Definition \ref{essential}). 

\noindent
(3)\,\,
For 
$\A':=\A \setminus \{H\}$, we say that $\A'$ is \textbf{NMPD} along $H$ 
if 
$\A'_p:=\{H \in \A' \mid 
p \in H\}$ is not maximal for all $p \in H$ with $\codim p \ge 3$.
\label{NMPD}
\end{define}

For example, 
$\A$ is NMPD along $H$ if $\pd \A$ is not maximal, and 
$\A$ is locally free (along $H$). 

To state a generalized addition-deletion theorems for projective dimensions, let us recall the Euler restriction.
For $H \in \A$, the 
\textbf{Euler restriction map} $\rho^H:D(\A) \rightarrow D(\A^H)$ is defined just by taking 
modulo $\alpha_H$. See Proposition \ref{ER} for details. Also, for the \textbf{localization} 
$$
\A_X:=\{H \in 
\A \mid X \subset H\}
$$
for $X\in L(\A^H)$, let $\rho_X^H:D(\A_X) 
\rightarrow 
D(\A^H_X)$ be the Euler restriction of $\A_X$. 
We say that $\rho^H$ is \textbf{locally surjective in codimension three along $H$} if $\rho_X^H$
is surjective for all $X \in L_2(\A^H)$.
Then we can give generalizations of Theorem \ref{adddel}
in the following manner:

\begin{theorem}[Addition theorem for projective dimensions]
Assume that $\rho^H$ is locally surjective in codimension three for $H \in \A$.
Let $\pd \A'=k$ and assume that $\A'$ is NMPD along $H$. 
\begin{itemize}
\item[(1)]
If $\pd \A^H = k$, then $\pd \A \le k+1$.
\item[(2)]
If $\pd \A^H > k$, then $\pd \A=\pd \A^H+1$.
\item[(3)]
If $\pd \A^H < k$, then $\pd \A=\pd \A'>\pd \A^H$.

\end{itemize}
In particular, the assumptions hold true if the $b_2$-equality holds for $(\A,H)$, and in that case 
the stronger statement (Theorem \ref{main}) holds.
\label{pdadd}
\end{theorem}

\begin{theorem}[Deletion theorem for projective dimensions]
Assume that $\rho^H$ is locally surjective in codimension three along $H \in \A$, 
$\pd \A=k$ and $\A'$ is NMPD along $H$. 
\begin{itemize}
\item[(1)]
If $\pd \A^H = k-1$, then $\pd \A' \le k$.
\item[(2)]
If $\pd \A^H > k-1$, then $\pd \A'=\pd \A^H$.
\item[(3)]
If $\pd \A^H < k-1$, then $\pd \A'=\pd \A>\pd \A^H+1$.

\end{itemize}

\label{pddel}
\end{theorem}

\begin{theorem}[Restriction theorem for projective dimensions]
Let $\pd \A'=k$ and $\A'$ is NMPD along $H$. 
Assume that $\rho^H$ is locally surjective in codimension three along $H \in \A$. 

\begin{itemize}
\item[(1)]
If $\pd \A < \pd \A'=k$, then $\pd \A^H=k$.
\item[(2)]
If $\pd \A'=\pd \A=k$, then $\pd \A^H \le k$.
\item[(3)]
If $\pd \A > \pd \A'=k$, then $\pd \A^H=\pd \A-1$.
\end{itemize}
In particular, the assumptions hold true if the $b_2$-equality holds for $(\A,H)$.
\label{pdrest}
\end{theorem}




\subsection{Surjectivity of restriction maps}

When we apply several results on freeness like Theorems \ref{adddel}, \ref{Ziegler2}, \ref{Y} and so on, the key 
is the surjectivity of not only the Euler but also 
the Ziegler restrictions. Let us explain briefly. The Ziegler restriction is the restriction 
that remembers the information on multiplicities. Namely, for $X \in \A^H$, define 
$$
m^H(X):=|\{L \in \A \setminus \{H\} \mid 
L \cap H=X\}|.
$$
Then $(\A^H,m^H)$ is called the \textbf{Ziegler restriction} of $\A$ 
onto $H$. We can define the logarithmic module $D(\A^H,m^H)$ in the same manner as for $\A$ 
(see Definition \ref{ma}). It is well-known that for 
$$
D_H(\A):=\{\theta \in 
D(\A) \mid \theta(\alpha_H)=0\},
$$
and the Euler derivation $\theta_E:=\sum_{i=1}^\ell x_i \partial_{x_i}$, 
it holds that $D(\A)=
S\theta_E\oplus D_H(\A)$ (see Lemma \ref{split}), and there is a map 
$$
\pi^H:=\rho^H|_{D_H(\A)}:D_H(\A) \rightarrow D(\A^H,m^H)
$$
called 
the \textbf{Ziegler restriction map}. 
Then Ziegler proved in \cite{Z} that $\pi^H$ is surjective if $\pd \A=0$, and conversely, 
$\pd \A=0$ if $\pd (\A^H,m^H):=
\pd_{S^H}
D(\A^H,m^H)=0$ and $\pi^H$ is surjective (Theorem \ref{Ziegler2}).
By using this converse implication, Yoshinaga proved 
in \cite{Y1} a freeness criterion in terms of the 
Ziegler restriction 
(Theorem \ref{Y}). Also, when we prove or apply Theorem \ref{adddel}, the surjectivity of $\rho^H$ is very important. Thus to determine when they are 
surjective is important, and Theorems \ref{pdadd}, \ref{pddel}, and \ref{pdrest} say that 
the local surjectivity at codimension three is important. 
The key to prove the above theorems are the following two results on the 
surjectivity. The first one is surprising in the sense that 
whether both $\pi^H$ and $\rho^H$ are surjective or not depends only on combinatorics if 
$\A$ is NMPD along $H$. 

\begin{theorem}

\begin{itemize}
\item[(1)]
Assume that $\A'$ is NMPD along $H$.
Then the Euler restriction map $\rho^H$ 
is 
surjective if and only if $\rho^H$ is locally surjective in codimension three along 
$H$.

\item[(2)]
Assume that $\A$ is NMPD along $H$.
Then 
the $b_2$-equality holds for 
$(\A,H)$ 
if and only if the Euler and the Ziegler restriction maps $\rho^H$ and 
$\pi^H$ onto $H$ are 
surjective.

\item[(3)]
Assume that $ \A'$ is NMPD along $H$.
Then 
the Euler restriction map $\rho^H$ onto $H$ is 
surjective if 
the $b_2$-equality holds for 
$(\A,H)$.
\end{itemize}
\label{surj}
\end{theorem}

Before the next step recall the multi-$b_2$-equality. We say that 
the \textbf{multi-$b_2$-equality} holds for $(\A,H)$ if 
$$
b_2^0(\A):=b_2(\A)
-|\A|+1=b_2(\A^H,m^H).
$$
See Proposition \ref{b2multi} on $b_2(\A^H,m^H)$. Then the result in \cite{AY} shows that 
$\A$ is free if and only if $(\A^H,m^H)$ is free and the multi-$b_2$-equality holds for $(\A,H)$ (see Theorem \ref{Y}). 
By using this equality, we can characterize 
the surjectivity of
$\pi^H$.

\begin{theorem}
Let $\pd \A=k \le \ell-3$ and $\A$ is NMPD along $H$.
Then the multi-$b_2$-equality holds for 
$(\A,H)$ if and only if the Ziegler restriction map $\pi^H$ onto $H$ is 
surjective.
\label{pisurj}
\end{theorem}

Let us show how to use Theorem \ref{pisurj} in the following example, that gives the case when the Ziegler restriction is surjective even though $\A$ is not free.

\begin{example}
Let us define $\B$ in $\R^4$ by 
$$
Q(\B):=x_1x_2x_3x_4(x_1-x_2)(x_2-x_3)(x_1-x_3)(x_1-x_4)(x_2-x_3-x_4)=0.
$$
First let us show that $\pd \B=1$. Let $L:x_2-x_4=0$ and 
$\A:=\B \cup \{L\}$. Then by using Theorem \ref{adddel} we can show that 
$\pd \A=0$ and $\exp(\A)=(1,3,3,3)$. Also, if $\pd \B=0$, then $\exp(\B)=(1,2,3,3)$. However, we can compute that $b_2(\B)=b_2(\A)-|\A^L|=30$. Hence this cannot occur. Now applying Theorem \ref{SPOG}, it holds that 
$\pd \B=1$. Hence Theorem \ref{NMPDcond} shows that $\B$ is NMPD along $H$. 

Next 
let $H:=\{x_4=0\} \in \B$. Then we can compute, by using Proposition \ref{b2multi},
$$
b_2^0(\B)=22=b_2(\B^H,m^H).
$$
Hence the multi-$b_2$-equality holds for $(\B,H)$.
Since $\ell-3=1 \ge 1=\pd \B$, Theorem \ref{pisurj} shows that the Ziegler restriction map 
$\pi^H:D_H(\B) \rightarrow D(\B^H,
m^H)$ is surjective without algebraic computations. This example shows that 
$\pi^H$ can be surjective even if $\B$ is not free. 
\end{example}

\subsection{Free case and the completion of the addition-deletion throrems}

When two of $\A,\A',\A^H$ are free, then we can apply Theorem \ref{adddel} 
with the $b_2$-equality. 
However, even from this viewpoint, Theorem \ref{adddel} is not complete. For example, 
it is shown in \cite{A5} that $\pd \A' \le 1$ if $\pd\A=0$ independent of 
$\pd \A^H$, but we do not know the converse, i.e., $\pd \A'=0$ and $\pd \A^H=0$ imply something on $\pd \A$ or not. This has been one of mysteries in free arrangement theory. For example, when $\ell=3$, then $\pd \A^H=0$ since the logarithmic derivation module is reflexive. Thus $\pd \A'=\pd \A^H=0$ do not imply $\pd \A=0$ in general. However, since it is reflexive again, $\pd \A\le 1$ in this case. 
Our first main result on these viewpoints is to complete this part for all $\ell \ge 3$. In fact, this 
result when $\ell=3$ is always true in the follwoing sense:

\begin{theorem}[$\pd$-version of the addition-deletion theorems]
Assume that two of $\pd \A,\ \pd \A',\ \pd \A^H$ are zero. Then 
the rest one is at most one. More precisely, 
\begin{itemize}
\item[(1)]
$\pd \A=\pd \A'=0$ implies that $\pd \A^H=0$. 
\item[(2)]
$\pd \A=\pd \A^H=0$ implies that $\pd \A'\le 1$. Explicitly, 
$\pd \A'=0$ if $\chi(\A^H;t) \mid \chi(\A;t)$, and 
$\pd \A'=1$ if $\chi(\A^H;t) \nmid \chi(\A;t)$. For the latter, $\A'$ is SPOG (Definition \ref{SPOGdef}). 
\item[(3)]
$\pd \A'=\pd \A^H=0$ implies that $\pd \A \le 1$. Explicitly, 
$\pd \A=0$ if $\chi(\A^H;t) \mid \chi(\A;t)$, and 
$\pd \A=1$ if $\chi(\A^H;t) \nmid \chi(\A;t)$. 
\end{itemize}
\label{main100}
\end{theorem}

Theorem \ref{main100} (1) follows immediately from Theorem \ref{adddel}, and 
(2) follows from Theorem \ref{SPOG}. Thus the main part is (3), which can be 
stated more details as follows:

\begin{theorem}
Let $H \in \A$, $\A':=\A \setminus \{H\}$. Assume that 
$\A'$ is free. Then $\A^H$ is free if and only if $\pd \A \le 1$. 
More precisely, 
\begin{itemize}
\item[(1)]
$\pd \A=0$ if both $\A'$ and $\A^H$ are free, and $\chi(\A^H;t) \mid \chi(\A;t)$. 
\item[(2)]
$\pd \A=1$ if both $\A'$ and $\A^H$ are free, and $\chi(\A^H;t) \nmid \chi(\A;t)$. 
\end{itemize}
\label{add1}
\end{theorem}

Theorem \ref{main100} completes the projective dimension table 
when two of the three projective 
dimensions are zero. To prove them, we again need 
the surjectivity of $\rho^H$, but that follows with no 
condition if $\A'$ is free as follows:

\begin{theorem}[Free surjection theorem]
Assume that $\pd \A'=0$. Then $\rho^H$ is surjective.
\label{freesurj}
\end{theorem}

Let us check above by an 
example.

\begin{example}
Let $\A$ be 
$$
xyzwu(x+w)(x+y+z+w)=0.
$$
Since $\chi(\A^H;t) \nmid \chi(\A;t)$, $\A':=\A \setminus \{x+y+z+w=0\}$ is free, and 
$(\A')^{\{x+y+z+w=0\}}$ is free, Theorem \ref{freesurj} shows that 
$\rho^H$ 
is surjective, and Theorem \ref{main100} shows that 
$\pd \A=1$. Since $\A'$ and $\A^H$ are divisionally free (see Theorem \ref{SFDF}), 
i.e., their freeness are combinaorially determined, 
we know that $\pd \A=1$ is also 
combinatorial.
\end{example}

As we can see above, we have a new class of arrangements whose projective dimension ($=1)$ can be determied just by combinatorics.

\begin{define}
$\PDC_1^\ell\ (\ell \le 2)$ are empty sets, 
$\PDC_1^3$ consists of non-empty arrangements $\A$ in $\K^3$ such that 
either $\chi_0(\A;t):=\chi(\A;t)/(t-1)$ is irredudible over $\Z$, or 
there is $H \in \A$ such that $\A \setminus \{H\} \in \DF_{3}$, and 
$\chi_0(\A;|\A^H|-1) \neq 0$. For 
$\ell \ge 4$, $\PDC_1^\ell$ consists of $\A$ such that 
there is $H \in \A$ such that $\A' \in \SF_\ell$, 
$\A^H \in \SF_{\ell-1}$ and $\chi(\A^H;t) \nmid \chi(\A;t)$. Let
$$
\PDC_1:=\cup_{\ell \ge 3} \PDC_1^\ell
$$
be the class of \textbf{inductively of projective dimension one}.
\label{PDC1}
\end{define}

Then by definition and Theorems \ref{main} and \ref{add1}, we have the following:

\begin{theorem}
$\pd \A=1$ if $\A \in \PDC_1$. Moreover, if $L(\A) \simeq L(\B)$ and 
$\A \in \PDC_1$, then $\pd \B=1$, i.e., $\pd \A=1$ is combinatorial.
\label{PDC1thm}
\end{theorem}


%
\subsection{Combinatorially determined projective dimensions}
By using our new methods to compute projective dimensions, we can construct new classes in which 
projective dimension is 
combinatorial. However, in Theorem \ref{pdadd}, we need the surjectivity of 
$\rho^H$ in codimension three whose 
combinatorial dependency 
is not yet known. Thus we give a family of arrangements in $\K^3$ in which 
$\rho^H$ is combinatorial

\begin{define}
Let $\A$ be an arrangement in $\K^3$ and 
let $H \in \A,\ \A':=\A \setminus \{H\}$. We say that 
a pair $(\A,H)$ is 
\textbf{a pair of combinatorially surjective $\rho^H$} if 
the Euler restriction $\rho^L:D(\B) 
\rightarrow D(\B^L)$ is surjective for any pair $(\B,L)$ of the arrangement 
$\B$ in $\K^3$ and $L \in \B$ such that there is a lattice isomorphism $f:L(\A) 
\rightarrow L(\B)$ with $f(H)=L$.

Let $\mbox{CS}_3$ be the set of all pairs $(\A,H)$ of an arrangement and a hyperplane that are combinatorially surjective $\rho^H$ in $\K^3$.
\label{CS}
\end{define}

Now 
let us introduce a new class of arrangements in which the projective dimension of 
their logarithmic derivation modules are combinatorial. 
For the divisionally and stair-free arrangements $\mathcal{DF}$ and $\mathcal{SF}$, see Definitions \ref{DF} 
and \ref{SF}. 

\begin{define}
Let $k\ge0$ 
and the set of \textbf{inductively of projective dimension 
$k$ in $\K^i$}, denoted by  
$\PDC^i_k$, be the set of arrangements in $\K^i$ defined as follows. First, 
$\PDC^i_0\ (i \le 2)$ 
coincides with all the arrangements in that ambient space, and 
$\PDC^i_{k>0}=\emptyset$ if $i=0,1,2$. 
$\PDC_0^3=\mathcal{DF}_3$ and 
let us use Definition 
\ref{PDC1} to define $\PDC_1^i$. 
Now for $i \ge 4$, define $\PDC_k^i$ as follows:
\begin{itemize}
\item[(1)]
$\PDC_0^i=\mathcal{SF}_i$. 
\item[(2)]
Let $\PDC_1^i[0]:=\PDC_1^i$ 
\item[(3)]
$\PDC_k^i[0]\ (k>1)$ consists of 
$\A$ with $H \in \A$ such that 
$\A \setminus \{H\} \in \PDC_{0}^i$ and 
$\A^H \in \PDC_{k-1}^{i-1}$. 
\item[(4)]
$\PDC_k^i[1]$ consists of 
$\A$ with $H \in \A$ such that 
$(\A_X^e,H) \in \mbox{CS}_3$ for any $X \in L_2(\A^H)$, 
$(\A')_X^e \in \PDC_{<k-1}^{k+1} $ for any $X \in L_k(\A^H)$ with $k>2$,
$\A \setminus \{H\} \in \PDC_{k-j}^i$ with $j >2$ and 
$\A^H \in \PDC_{k-1}^{i-1}$. 
\item[(5)]
$\PDC_{\ell-2}^i[2]$ consists of 
$\A$ with $H \in \A$ such that the $b_2$-equality holds for $(\A,H)$, 
$(\A')_X^e \in \PDC_{<k-1}^{k+1} $ for any $X \in L_k(\A^H)$ with $k>2$, and 
$\A \setminus \{H\} \in \PDC_{\ell-2}^i[0]
\cup \PDC_{\ell-2}^i[1]\cup \PDC_{\ell-2}^i[2]$. 
$\PDC_{k}^i[2]=\emptyset$ if $k < \ell-2$. 

\item[(7)]
Now 
$$
\PDC_k^i:=
\PDC_k^i[0] \cup
\PDC_k^i[1] \cup
\PDC_k^i[2] .
$$
\end{itemize}
Finally, define the set of 
\textbf{inductively of projective dimension 
$k$} by 
$$
\PDC_k:=\cup_{i=0}^\infty \PDC_k^i.
$$
\label{PDC}
\end{define}

Actually we can show that in $\PDC$ projective 
dimension is combinatorial.

\begin{theorem}

\begin{itemize}
\item[(1)]
Let $\A \in \PDC_k^i$. Then $\pd \A=k$.

\item[(2)]
Let $\A \in \PDC_k^i$ and $\B$ be another arrangement. If 
$L(\A) \simeq L(\B)$, then $\pd \B=k$.
\end{itemize}
\label{PDCcombin}
\end{theorem}

Let us check it by an 
example.

\begin{example}
(1)\,\,
Let $\A'$ be 
$$
xyzw(x-y)(x-z)(x-w)(y-z)(y-w)(z-w)=0.
$$
This is free 
with $\exp(\A')=(1,2,3,4)$. It is famous that 
$\A' \in \DF_3 = \mathcal{SF}_3=\PDC_0^3$. 
Let $H:x+y+z+w=0$ and 
let 
$\A:=\A' \cup \{H\}$. Since 
$$
\chi_0(\A^H)=t^2-9t+26,
$$it holds that $\A^H \in \PDC_1^3$ by definition. Thus $\pd \A^H=1$. 
Hence Definition \ref{PDC} (2) shows that $\A \in \PDC_2^4$.
So Theorem 
\ref{PDCcombin} implies that, for any $\B$ with $L(\B) \simeq L(\A)$, it holds that 
$\pd \B=2$.

(2)\,\,
Let $\A'$ be 
$$
xyzwu(x+w)(x+y+z+w)=0.
$$
Since $\chi(\A^H;t) \nmid \chi(\A;t)$, $\B:=\A' \setminus \{x+y+z+w=0\} \in \mathcal{DF}_5$ and 
$(\A')^{\{x+y+z+w=0\}} \in \mathcal{DF}_4$, Theorem \ref{add1} shows that 
$\A' \in \PDC_1^5$. Let $H:y-w-u=0$ and $\A:=\A' \cup \{H\}$. Then 
$$
Q(\A^H,m^H)=xyzw(x+w)(x+y+z+w)(y-u).
$$
Thus applying Theorem \ref{CS} to $(\A_X,H)$
for all $X \in L_2(\A^H)$, we can show that 
$(\A_X,H) \in \mbox{CS}_3$, and 
$\A' \in \PDC_1^4$, thus $\A'$ is NMPD along $H$ by Proposition \ref{NMPDcond}. 
Also, $\A^H \in \mathcal{PDC}_2^4$. Hence 
Definition \ref{PDC} and Theorem \ref{PDCcombin} show that 
$\A \in \PDC_3^5$, thus for any $\C$ with $L(\C) \simeq L(\A)$, it holds that 
$\pd \C=3$.
\end{example}

The organization of this article is as follows. In \S2 we collect several definitions and 
results for the proof of main results that is done in 
\S3. In \S4 we prove Theorem {PDCcombin} and some conditions for $\rho^H$ to be 
surjective. \S5 is devoted for several variants and applications of our main results. In \S6 we pose 
several problems related to results in this article.
\medskip

\noindent
\textbf{Acknowledgements}. The author is
partially supported by JSPS KAKENHI 
Grant-in-Aid for Scientific Research (B) Grant Number JP16H03924. 
The author is grateful to M. DiPasquale for the discussion in \cite{D} on 
Example \ref{noncombex}.

\section{Preliminaries}
In this section we recall several definitions and results. 
Let $\A$ be a central \textbf{arrangement (of hyperplanes)} in 
$V=\K^\ell$, i.e., a finite set of linear hyperplanes in $V$. 
Let 
$S:=\mbox{Sym}^*(V^*) \simeq \K[x_1,\ldots,x_\ell]$ be the 
coordinate ring of $V$, and $\Der S =\oplus_{i=1}^\ell S \partial_{x_i}$ the 
module of $\K$-linear $S$-derivations. For a fixed defining linear form $\alpha_L \in V^*$ of 
$L \in \A$, the logarithmic derivation module $D(\A)$ is defined by 
$$
D(\A):=\{\theta \in \A \mid 
\theta(\alpha_H) \in S \alpha_H\ (\forall H \in \A)\}.
$$
$D(\A)$ is a reflexive $S$-graded module, and not free in general. Hence we say that $\A$ is 
\textbf{free} with \textbf{exponents} 
$\exp(\A)=(d_1,\ldots,d_\ell)$ if $D(\A)$ is a free $S$-module with 
a homogeneous basis $\theta_1,\ldots,\theta_\ell$ of degree $\deg \theta_i=d_i$ 
for $i=1,\ldots,\ell$. Here $\deg \theta_i:=\deg \theta_i(\alpha)$ for some linear form 
$\alpha$ with $\theta_i(\alpha) \neq 0$. Since $D(\A)$ is reflexive, the following is well-known:

\begin{lemma}
$$
0 \le \pd_S D(\A) \le \ell-2.
$$
\label{pd-2}
\end{lemma}

When $\A$ is not empty, the degree one element $\theta_1 \in D(\A)$ can be chosen from 
a set of generators for $D(\A)$, called the \textbf{Euler 
derivation} $\theta_E=\sum_{i=1}^\ell x_i \partial_{x_i}$ which is always contained in $D(\A)$. 
Also, for $H \in \A$, define 
$D_H(\A):=\{\theta \in D(\A) \mid \theta(\alpha_H)=0\}$. Then we have the following:

\begin{lemma}[Lemma 1.33 in \cite{Y3} for example]
$$
D(\A)=S\theta_E \oplus D_H(\A),
$$
\label{split}
\end{lemma}

Thus we have the following:

\begin{lemma}
It  holds that 
$\pd \A=\pd_S D_H(\A)$.
\label{decomp}
\end{lemma}

For a direct sum decomposition $V=V_1 \oplus V_2$ of the vector space $V$, let 
$\A_i$ be an arrangement in $V_i$ and let $S_i$ be the coordinate ring of $V_i\ (i=1,2)$. 
Then $S_1 \otimes_\K S_2=S$. Define $$
\A_1 \times \A_2:=\{H_1 \oplus V_2 \mid H_1 \in \A_1\} \cup \{V_1 \oplus H_2 \mid H_2 \in \A_2\}.
$$
Then we have the following.

\begin{prop}[\cite{OT}, Proposition 4.14]

$$
D(\A_1 \times \A_2)=S\cdot D(\A_1) 
\oplus S \cdot D(\A_2).
$$
\label{directsum}

\end{prop}

To state the advantage of free arrangements, let us introduce combinatorics and topology of arrangements. 

\begin{define}
(1)\,\,The \textbf{intersection lattice} $L(\A)$ of $\A$ is defined by 
$$
L(\A):=\{ \cap_{H \in \B} H \mid \B  \subset \A\}. 
$$ 
A partial order in $L(\A)$ is equipped with the 
reverse inclusion. Let $L_k(\A):=\{X \in L(\A)\mid \codim_V X=k\}$, and 
\begin{eqnarray*}
L_{\ge k}(\A):&=&\{ X \in L(\A) \mid \codim X \ge k\},\\
L_{\le k}(\A):&=&\{ X \in L(\A) \mid \codim X \le k\}.
\end{eqnarray*}

(2)\,\,
The \textbf{M\"{o}bius function} $\mu : L(\A) 
\rightarrow \Z$ is defined by $\mu(V)=1$, and by 
$\mu(X):=-\sum_{X  \subsetneq Y \subset V,\ Y \in L(\A)} \mu(Y)$ for $X \in L(\A) \setminus \{V\}$. 
\label{IL}
\end{define}

\begin{define}
The \textbf{characteristic polynomial} $\chi(\A;t)$ of $\A$ is defined by 
$$
\chi(\A;t):=\sum_{X \in L(\A)} \mu(X)t^{\dim X}=:\sum_{i=0}^\ell
(-1)^{i}b_i(\A)t^{\ell-i}.
$$
and 
the \textbf{Poincar\`{e} polynomial} $\pi(\A;t)$ of $\A$ is defined by 
$$
\pi(\A;t):=\sum_{X \in L(\A)} \mu(X) (-t)^{\codim X}
=\sum_{i=0}^\ell b_i(\A)t^{i}.
$$
It is known that $\pi(\A;t)=\mbox{Poin}(V \setminus 
\cup_{L \in \A} L;t)$ when $\K=\CC$. 

When $\A \neq \emptyset$, it is known that $\chi(\A;t)$ is divisible by 
$t-1$. Define 
$$
\chi_0(\A;t):=\chi(\A;t)/(t-1)=\sum_{i=0}^{\ell-1}
(-1)^ib_i^0(\A)t^{\ell-1-i}=\sum_{i=0}^{\ell-1}
(-1)^ib_i(d\A)t^{\ell-1-i},
$$
where $d\A$ is the deconing of $\A$ by any line $H \in \A$.
\label{char}
\end{define}

We can define two fundamental operations to create new arrangements from a given $\A$ and $X \in L(\A)$. 

\begin{define}
For $X \in L(\A)$,
define 
\begin{eqnarray*}
\A_X:&=&\{H\in \A \mid X \subset H\},\\
\A^X:&=&\{L \cap X \mid L \in \A \setminus \A_X\}.
\end{eqnarray*}
$\A_X$ is called the \textbf{localization}, and $\A^X$ is the \textbf{restriction} of $\A$ onto $X$. The former is an arrangement in $V$, and the latter in $X \simeq \K^{\dim X}$.

More 
generally, for a homogeneous prime ideal $p \subset S$, let 
$$
\A_p:=\{H \in \A \mid p \in H\}
$$
be the localization at the point $p$. 
\label{localrest}
\end{define}

Localizations play important roles in this article. To use localizations, 
it is important to see the essential part of them as follows:

\begin{define}
For $T:=\cap_{H \in \A} H \neq 0$, $\A$ can be expressed as a direct product of some 
arrangement $\B$ and $\emptyset_k$ as $\A=\B \times \emptyset_k$, where 
$k:=\dim T$, $\emptyset_k$ is the \textbf{empty 
arragngement} in $\K^k$, and $\B$ is an arrangement in $V/T$ with $\cap_{L \in \B} L=0$. $\B$ is called the \textbf{essentialization} 
of $\A$, and denoted by $\B=\A^e$ in this article.
\label{essential}
\end{define}

We use essentialization frequently when we consider the localizaion $\A_X$ of $\A$ at $X \in L_k(\A)$, i.e., in 
this case $\A_X=\A_X^e \times \emptyset_{\ell-k}$. Here $\A^e_X$ is an arrangement in 
$V/X\simeq \K^{k}$. By the definition of the localization and the projective resolution, the following is clear (see also Corollary 5.2 in \cite{To}).

\begin{lemma}
Let $p \subset S$ be a homogeneous prime ideal. Then 
$$
\pd \A_p \le \pd \A.
$$ 
Thus 
for  $X \in \A$, it holds that 
$$
\pd \A_X \le \pd \A.
$$
In particular, $\A_p$ and $\A_X$ are free if $\A$ is free.
\label{pdlow}
\end{lemma}

Now let us reformulate the first restriction that plays a key role in this article.

\begin{prop}[Euler restriction]
There is an exact sequene 
$$
0 \rightarrow 
D(\A \setminus \{H\}) 
\stackrel{\cdot \alpha_H}{\rightarrow }
D(\A)
\stackrel{\rho^H}{\rightarrow } D(\A^H).
$$
Here $\rho^H(\theta)(\overline{f}):=\overline{\theta(f)}$ modulo $\alpha_H$ for the image $\overline{f} 
\in S^H$ of $f \in S$ in $S/\alpha_H S$ is 
called the 
\textbf{Euler restriction}.
\label{ER}
\end{prop}

We may relate the exponents of free arrangements, combinatorics and 
topology as follows:

\begin{theorem}[Terao's factorization, \cite{T2}]
Assume that $\A$ is free with $\exp(\A)=(d_1,\ldots,d_\ell)$. Then 
$\chi(\A;t)=\prod_{i=1}^\ell(t-d_i)$.
\label{Teraofactorization}
\end{theorem}

For the analysis of the freeness, the following is of the most importance.

\begin{theorem}[Terao's polynomial $B$-theory, \cite{T1}]
Let $H \in \A$, 
$\A':=\A 
\setminus \{H\}$ and let us define the homogeneos degree 
$(|\A'|-|\A^H|)$-polynomial $B$ by 
$$
B:=\prod_{X \in \A^H} \alpha_{\nu(X)},
$$
where $\nu:\A^H \rightarrow \A$ is a section such that 
$\nu(X) \cap H=X$. The polynomial $B$ is called \textbf{Terao's polynomial $B$}. Then for an arbitrary $\theta \in D(\A')$, it holds that 
$$
\theta(\alpha_H) \in (\alpha_H,B).
$$
Thus, $\theta \in D(\A')$ is in $D(\A)$ if $\deg \theta <|\A'|-|\A^H|$.
Moreover, if there is $\varphi \in 
D(\A')$ such that $\deg \varphi=|\A'|-|\A^H|$ and that 
$\varphi \not \in D(\A)$, then for $\theta \in D(\A')$, there is 
$f \in S$ such that $\theta -f \varphi \in D(\A)$. Thus 
$$
D(\A)=D(\A')+S\cdot \varphi.
$$
\label{B}
\end{theorem}

Next let us introduce the theory for multiarrangements, which was introduced by Ziegler 
in \cite{Z}.

\begin{define}[\cite{Z}]

(1)\,\,
A pair $(\A,m)$ is a \textbf{multiarrangement} if 
$m:\A \rightarrow \Z_{\ge 1}$. For $X \in L(\A)$, let $m_X:=m|_{\A_X}$ and 
the pair $(\A_X,m_X)$ is called the \textbf{localization} of 
$(\A,m)$ at $X$. Let 
$$
Q(\A,m):=\prod_{H \in \A} \alpha_H^{m(H)}.
$$

\noindent
(2)\,\, 
The \textbf{logarithmic derivation module} $D(\A,m)$ of $(\A,m)$
is defined by 
$$
D(\A,m):=\{\theta \in \Der S \mid \theta(\alpha_H) \in S \alpha_H^{m(H)}\ (\forall H \in \A)\}.
$$
\label{ma}
\end{define}

We may define its \textbf{freeness} and \textbf{exponents} in the same way as 
for $m \equiv 1$. Also, let 
$$
\pd (\A,m):=\pd_S D(\A,m).
$$
Note that if $\ell=2$ then $(\A,m)$ is free since $D(\A,m)$ is always reflexive. 
The most fundamental
criterion for the freeness is the following:

\begin{theorem}[Saito's criterion, \cite{Sa}, \cite{Z}]
Let $\theta_1,\ldots,\theta_\ell \in D(\A,m)$ 
and let $M=(\theta_i(x_j))$ be the $(\ell \times \ell)$-matrix. Then 
$$
\det M \in S Q(\A,m).
$$
Moreover, $\A$ is free with 
basis $\theta_1,\ldots,\theta_\ell$ if and only if 
$$
\det M=Q(\A,m)
$$
up to non-zero scalar. 
\label{saito}
\end{theorem}

Also, we can define the characteristic polynomial $\chi(\A;t)=\sum_{i=0}^\ell b_i(\A,m)t^{\ell-i}$ of $(\A,m)$ in a algebraic way, see 
\cite{ATW} for details. In general it is very difficult to compute $b_i(\A,m)$ except for $b_1(\A,m)=|m|:=\sum_{H \in \A}
m(H)$. However, for the second Betti number $b_2(\A,m)$ of $(\A,m)$, we have the following way to compute.

\begin{prop}[\cite{ATW}, 
Corollary 4.4]
\begin{itemize}
\item[(1)]
Assume that $\ell=2$ and $\exp(\A,m)=(d_1,d_2)$. Then $b_2(\A,m)=d_1d_2$.

\item[(2)]
For $X\in L_2(\A)$, let $\exp(\A_X,m_X)=(d_1^X,d_2^X,0,\ldots,0)$. 
Then 
$$
b_2(\A,m)=\sum_{X \in L_2(\A)} b_2(\A_X,m_X)=\sum_{X 
\in
L_2(\A)} d_1^X d_2^X.
$$
\end{itemize}

\label{b2multi}
\end{prop}

Also, the following local functoriality of $D(\A,m)$ is important too.

\begin{lemma}[\cite{OT}]
Let $p$ be a homogeneous prime ideal of $S$. 
Then 
$$
D(\A,m)_p =D(\A_p,m|_{\A_p})_p.
$$
In particular, for a generic point $p$ of $X \in L(\A)$, it holds that 
$$
D(\A,m)_p =D(\A_X,m_X)_p.
$$
\label{localfunctor}
\end{lemma}

We can construct the multiarrangement canonically from an arrangement $\A$.

\begin{define}
For an arrangement $\A$ in $\K^\ell$ and $H \in \A$, define 
$m^H(X):=|L \in \A \setminus \{H\} \mid 
L \cap H=X\}|$ for $X \in \A^H$. The pair $(\A^H,m^H)$ is 
called the \textbf{Ziegler restriction} of $\A$ onto $H$. Also, there is the 
\textbf{Ziegler restriction map} 
$$
\pi^H=\pi:D_H(\A) \rightarrow D(\A^H,m^H)
$$
by taking modulo $\alpha_H$. 
Equivalently, 
$$
\pi^H:=\rho^H|_{D_H(\A)}.
$$
\label{Zieglerrest}
\end{define}

Whether the Euler and Ziegler restrictions are surjective or not is 
difficult to see. 

\begin{example}
\label{ex20}
(1)\,\,
Let $\A$ be defined by $xyz(x+y+z)=0$, whose $D(\A)$ has a generator of degrees 
$1,2,2,2$. Then for $H:=\ker z$, $\A^H$ is defined by $xy(x+y)=0$, whose exponents 
are $(1,2)$. In this case, it is easy to show that $\rho^H$ is surjective. However, $\pi^H$ is not 
surjective since by 
Yoshinaga's criterion (see Theorem \ref{Y} below), 
$$
\codim \coker \pi^H=b_2^0(\A)-1\cdot 2=3-2=1.
$$

(2)\,\,
Let $\A$ be a $3$-arrangement obtained as the coning of the affine arrangement consisting of 
edges and diagonals of the regular pentagon. It is known to be free with exponents $(1,5,5)$, and 
$|\A^H|=5$ for all $H \in \A$. 
In this case, $\pi^H$ is surjective by the freeness of $\A$ and Theorem \ref{Ziegler2}. However, $\rho^H$ is not since $\exp(\A^H)=(1,4)$. 
\end{example}

A remarkable property of the Ziegler restriction map is the following.

\begin{theorem}[\cite{Z}]
(1)\,\,
Assume that $\A$ is free with $\exp(\A)=(1,d_2,\ldots,d_\ell)$. Then 
for any $H \in \A$, the Ziegler restriction $(\A^H,m^H)$ is also free with 
$\exp(\A^H,m^H)=(d_2,\ldots,d_\ell)$. In particular, $\pi^H$ is surjective.

\noindent
(2)\,\,
Conversely, $\A$ is free if and only if $\pi_H$ is 
surjective, and  $(\A^H,m^H)$ is free.
\label{Ziegler2}
\end{theorem}

Moreover, a converse of Theorem \ref{Ziegler2} holds true with 
additional conditions.

\begin{theorem}[Yoshinaga's criterion, \cite{Y1}, \cite{Y2}, \cite{AY}]

In the notation of Definition \ref{Zieglerrest}, 
the \textbf{multi-$b_2$-inequality} 
$$
b^0_2(\A) \ge b_2(\A^H,m^H).
$$
holds. 
Moreover, $\A$ is free if and only if the above inequality is the equality 
(\textbf{multi-$b_2$-equality}), and 
$(\A^H,m^H)$ is free. In particular, $b_2^0(\A)=b_2(\A^H,m^H)$ if and only if  $\pi^H$ is locally 
surjective in codimension three along $H$. Also, when $\ell=3$, it holds that 
$$
\A\ \mbox{is free} \iff 
b^0_2(\A) = b_2(\A^H,m^H) \iff \pi^H \ \mbox{is surjective}.
$$
\label{Y}
\end{theorem}

%

The following is a generalization of Terao's addition-deletion theorem in terms of 
the $b_2$-equality.

\begin{theorem}[$b_2$-inequality and the division theorem, \cite{A2}]
Let $H \in \A$. 
\begin{itemize}
\item[(1)]
It holds that 
$$
b_2(\A) \ge b_2(\A^H)+|\A^H|(|\A|-|\A^H|),
$$
and the equality (called the \textbf{$b_2$-equality}) implies that $\A_X:=\{H \in \A \mid H \supset X\}$ 
is free for all $X \in L(\A^H)$ with $\codim_V X=3$. 

\item[(2)]The $b_2$-equality implies the multi-$b_2$-equality.
\item[(3)] 
$\pd \A=0$ if for some $H \in \A$,
\begin{itemize}
\item[(i)]
$\A^H$ is free, and 
\item[(ii)]
the $b_2$-equality holds for $(\A,H)$.
\end{itemize}
\item[(4)]
Assume that the multi-$b_2$-equality holds for $(\A,H)$. Then $\pd (\A^H,m^H)=0$ if 
$\pd \A^H=0$ .
\end{itemize}
\label{division}
\end{theorem}

Let us divide these equalities as follows:

\begin{define}
Let $H \in \A$. Then we say that the 
\textbf{upper $b_2$-equality}, or the \textbf{multi-$b_2$-equality} holds for $(\A,H)$ if 
$$
b_2^0(\A)=b_2(\A^H,m^H).
$$
Also, 
we say that the 
\textbf{lower $b_2$-equality} holds for $(\A,H)$ if 
$$
b_2(\A^H,m^H)=b_2(\A^H)+(|\A^H|-1)(|\A|-|\A^H|-1).
$$
\label{ulb2}
\end{define}

The relation among these
equalities are as follows:

\begin{prop}
(1)\,\,
In general, it holds that 
$$
b_2^0(\A)\ge b_2(\A^H,m^H).
$$
and that 
$$
b_2(\A^H,m^H) \ge b_2(\A^H)+(|\A^H|-1)(|\A|-|\A^H|-1).
$$

\noindent
(2)\,\,
The $b_2$-equality 
$$
b_2(\A)=b_2(\A^H)+
|\A^H|(|\A|-|\A^H|)
$$
is equivalent to 
$$
b_2^0(\A)=b_2(\A)-|\A|+1=
b_2(\A^H)+
(|\A^H|-1)(|\A|-|\A^H|-1).
$$

\noindent
(3)\,\,
The $b_2$-
equality holds if and only if both the upper and 
lower $b_2$-equalities hold.
\label{b2prop}
\end{prop}

The (multi) $b_2$-equalities defined above are local in the following sense.

\begin{prop}
(1)\,\,
Assume that the (multi-)$b_2$-equality holds for $(\A,H)$. Then the same holds for 
$(\A_X,H)$, where $X \in L(\A^H)$.

\noindent
(2)\,\,
The (multi-)$b_2$-equality holds for $(\A,H)$ if and only if the same holds for 
$(\A_X,H)$, where $X$ runs all the element in $L_2(\A^H)$.

\label{b2local}
\end{prop}

\noindent
\textbf{Proof}. 
(1)\,\,
Since the proof is the same, we prove when 
$$
b_2^0(\A)=b_2(\A^H,m^H)
$$
holds. Let $X \in L(\A^H)$. Then by Proposition \ref{b2multi}, it holds that 
$$
b_2^0(\A_X)=\sum_{dX \subset Y \in L(d\A)} b_2(\A_Y)
$$
and 
$$
b_2(\A^H_X,m^H_X)=\sum_{X \subset Y \in L_2(\A^H)} b_2(\A_Y^H,m^H_Y).
$$
Here $d\A$ and $dX$ are the deconing of $\A$ and $X$ with respect to $H$, i.e., 
$d\A=\A|_{\alpha_H=1}$ and $dX=X|_{\alpha_H=1}$. 
Note that 
\begin{eqnarray*}
b_2^0(\A)&=&
\sum_{dX \subset Y \in L_2(d\A)} b_2(\A_Y) +
\sum_{dX \not \subset Y \in L_2(d\A)} b_2(\A_Y)\\
=b_2(\A^H,m^H)&=&\sum_{X \subset Y \in L_2(\A^H)} b_2(\A_Y^H,m^H_Y)+
\sum_{X \not \subset Y \in L_2(\A^H)} b_2(\A_Y^H,m^H_Y).
\end{eqnarray*}
Since
$$
b_2(d\A_Y)=b_2^0(\A_Y) \ge b_2(\A_Y^H,m_Y^H)
$$
by Theorem \ref{Y}, this inequality is the equality, and Proposition \ref{b2multi} 
implies that 
$b_2^0(\A_X)=b_2(\A_X^H,m^H_X)$. 

(2)\,\, Apply (1) and Proposition \ref{b2multi}. 
\owari
\medskip

\begin{prop}[cf. \cite{AY0}]

For a coherent sheaf $E$ on $\P^{\ell-1}$, let 
$\Gamma_*(E):=\oplus_{d\in \Z} H^0(E(d))$. Then 
\begin{itemize}
\item[(1)]
$\Gamma_*(\widetilde{D(\A)})=D(\A)$,
\item[(2)]
$\Gamma_*(\widetilde{D_H(\A)})=D_H(\A)$, and 
\item[(3)]
$\Gamma_*(\widetilde{D(\A,m)})=D(\A,m)$.
\end{itemize}
\label{globalsection}
\end{prop}

\noindent
\textbf{Proof}. We first prove (3), then (1) and (2) follow immediately by the definition. 
Let $E:=\widetilde{D(\A,m)}$ be a sheaf on $\P^{\ell-1}$. Let $K$ be the quotient field of $S$. Then 
$E \subset K^{\ell}$ as a sheaf, thus $\theta \in \Gamma_*(E) \in K^\ell$. Thus $\theta$ can be expressed as 
$$
\theta=\displaystyle \frac{\theta_1}{x_1^{d_1}}=
\cdots=\displaystyle \frac{\theta_\ell}{x_\ell^{d_\ell}} \in \K^\ell
$$
for $\theta_i \in D(\A,m)$. Thus $\theta \in \Der S$. Let $H \in \A$. Then we may assume that $(\alpha_H,x_1)=1$. Thus 
$\theta_1(\alpha_H) =\exists f \alpha_H^{m(H)}$ implies that $\theta \in D(\A,m)$. Since $D(\A,m) \subset \Gamma_*(E)$ is well-known, 
we complete the 
proof.\owari
\medskip

\begin{theorem}[\cite{A4}, Theorem 1.6 and Corollary 1.7]
Let $\A$ be an $\ell$-arrangement and 
$(\A^H,m^H)$ the Ziegler restriction of $\A$ onto $H \in \A$. Assume that 
the $b_2$-equatity
$$
b_2(\A)=b_2(\A^H)+|\A^H|(|\A|-|\A^H|)
$$
holds true. Let 
$\theta_E^H:=(Q(\A^H,m^H)/Q(\A^H))\theta_E \in D(\A^H,m^H)$, and 
let $\pi^H:D_H(\A) \rightarrow D(\A^H,m^H)$ be the Ziegler 
restriction.

(1)\,\,
Then there are generators $\theta_E,\theta_1,\ldots,\theta_s$ for 
$D(\A^H)$ such that $\theta_E^H,
\theta_1,\ldots,\theta_s$ form a generator for $D(\A^H,m^H)$. 

(2)\,\,
Assume that $\A$ 
is free, and $\A^H$ is not free. Then there are generators $
\theta_E,\theta_2,\ldots,\theta_\ell$ for $D(\A^H)$ such that 
the preimage of $\theta_2,\ldots,\theta_\ell$ in $D_H(\A)$ by $\pi^H$ 
form a free basis for
$D_H(\A)$, and the relation among them is in the degree 
$|\A|-|\A^H|$ of the form 
$
\theta_E^H=\sum_{i=2}^\ell f_i \theta_i,
$
and no other relation exists.
\label{b2gen}
\end{theorem}

Also, the following relation between Betti numbers and Chern classe are important.

\begin{prop}[\cite{DS}, Proposition 5.18]
$b_i^0(\A)=c_i(\widetilde{D_0(\A)})$ for $i=0,1,2$.
\label{DS}
\end{prop}

\begin{prop}
Assume that there is the $b_2$-equatity for $(\A,H)$. 
Then $\rho^H$ is surjective if $\pi^H$ is.
\label{pirho}
\end{prop}

\noindent
\textbf{Proof}. 
By Theorem \ref{b2gen}, there is a minimal generator 
$\theta_E,\theta_1,\ldots,\theta_s$ for $D(\A^H)$ such that 
$\theta_E^H,\theta_1,\ldots,\theta_s$ form a generator for 
$D(\A^H,m^H)$. 
Since $\pi^H$ is surjective, there are derivations 
$\varphi_i\ (i=1,\ldots, s)$ in $D_H(\A)$ such that $\pi^H(\varphi_i)=\rho^H(\varphi_i)=\theta_i$. Since $\theta_E \in D(\A)$ and 
$\rho^H(\theta_E) \in D(\A^H)$, it holds that 
$\rho^H$ is surjective. \owari
\medskip


The next result is used to determine a set of generators for $D_H(\A)$. 

\begin{theorem}[\cite{AD}, Theorem 5.1]
Let $H \in \A$ and 
$\pi^H:D_H(\A) \rightarrow D(\A^H)$ be the Ziegler restriction. 
If $\mbox{Im}(\pi^H)$ is generated by 
$\pi^H(\theta_2),\ldots,\pi^H(\theta_s)$ for 
$\theta_i \in D_H(\A)$, then $D(\A)$ is generated by 
$$
\theta_E,\theta_2,\ldots,\theta_s.
$$
\label{gen}
\end{theorem}

For the projective dimensions of free minus-one arrangements, we have the following 
explicit results.

\begin{theorem}[Theorem 1.4, \cite{A5}]
Let $\A$ be free, and $H \in \A$. If $\A':=\A \setminus \{H\}$ is not free, then 
$D(\A')$ has a minimal free resolution of the following form:
$$
0 \rightarrow 
S[-e_H-1] \rightarrow D(\A) \oplus S[-e_H] \rightarrow D(\A') \rightarrow 0.
$$
Here $e_H:=|\A|-|\A^H|-1$. 
In particular, $\pd \A' \le 1$ for all free arrangement $\A$ and $H \in \A$.
\label{SPOG}
\end{theorem}

\begin{define}[\cite{A5}]
We say that $\A$ is \textbf{plus-one generated (POG)} with 
$\POexp(\A)=(d_1,\ldots,d_\ell)$ and level $d$ if 
$D(\A)$ has a following minimal free resolution:
$$
0 \rightarrow 
S[-d-1] \stackrel{(f_1,\ldots,f_\ell,\alpha)}{\longrightarrow} \oplus_{i=1}^\ell S[-d_i] \oplus S[-d] \rightarrow D(\A) \rightarrow 0.
$$
We say that $\A$ is \textbf{strictly plus-one generated (SPOG)} with 
$\POexp(\A)=(d_1,\ldots,d_\ell)$ and level $d$ if 
$D(\A)$ is POG and $\alpha \neq 0$ in the above notation. Such a set of 
minimal homogeneous generators $\theta_1,\ldots,\theta_\ell,\varphi$ with $\deg 
\theta_i=d_i,\ \deg \varphi=d$ is called a \textbf{SPOG-generator}. 
\label{SPOGdef}
\end{define}

Let us recall two classes 
of free
arrangements in which the freeness depends only on $L(\A)$.

\begin{define}[\cite{A2}]
The class $\mathcal{DF}_\ell$ of arrangements in $\K^\ell$ consists of arrangements 
$\A$ such that, there is $X_i \in L_i(\A)\ (i=0,\ldots,\ell-2)$ such that 
$$
X_0=V \supset X_1 \supset \cdots \supset X_{\ell-2},
$$
and $b_2(\A^{X_{i-1}})=b_2(\A^{X_{i}})+|\A^{X_{i}}|(
|\A^{X_{i-1}}|-|\A^{X_{i}}|)$ for $i=1,\ldots,\ell-2$. The set 
$$
\mathcal{DF}:=\cup_{\ell \ge 0} \mathcal{DF}_\ell
$$
is called the set of \textbf{divisionally free arrangements}, and the above flag $\{X_i\}_{i=0}^{\ell-2}$ of 
$\A$ is called a \textbf{divisional flag} of $\A$.
\label{DF}
\end{define}

\begin{define}[Stair-free arrangements, \cite{A6}]
We say that $\A subset \B$ is connected by a \textbf{free path} if 
there is an order $\B \setminus \A=\{H_1,\ldots,H_s\}$ such that 
$\A \cup\{H_1,\ldots,H_i\}$ is free for all 
$1 \le i \le s$. We say that $\A$ and $\A^H\ (H \in \A)$ is \textbf{divisionally connected} 
if there is the $b_2$-equality for $(\A,H)$. 
The set $\mathcal{SF}_\ell$ consists of hyperplane arrangements $\A$ 
in $\K^\ell$ such that $\A$ is connected to $\phi$ by the divisional connected and free paths. 
$\A \in \mathcal{SF}$ is called a \textbf{stair-free arrangement of hyperplanes}. 
\label{SF}
\end{define}

For the above two classes of arrangements, the combinatorial freeness was shown as follows:

\begin{theorem}[\cite{A2}, \cite{A6}]
$\mathcal{SF} \supset \DF$, $\A \in \mathcal{SF}$ is free, and 
the Terao's 
conjecture is true in $\mathcal{SF}$.
\label{SFDF}
\end{theorem}

\section{Proofs}
For the Euler and Ziegler restrictions $\rho^H,\pi^H$ and $X \in L(\A^H)$, let 
$\rho^H_X$ be the Euler restriction $D(\A_X) \rightarrow D(\A^H_X)$ and 
$\pi_X^H$ be the Ziegler restriction $D_H(\A_X) \rightarrow D(\A^H_X,m_X^H)$. By Lemma \ref{localfunctor}, it holds that 
$$
(\rho^H)_p=(\rho_X^H)_p,\ (\pi^H)_p=(\pi_X^H)_p
$$
for a generic point $p \in X$.

To prove main results, let us introduce several results. First results are common facts on 
commutative algebra.

\begin{lemma}
Assume that the Euler restriction $\rho^H:D(\A) \rightarrow 
D(\A^H)$ is surjective. Then the localization 
$\rho_X^H:D(\A_X) \rightarrow D(\A_X^H)$ is surjective for 
$X \in L(\A^H)$.
\label{surjloc}
\end{lemma}

\noindent
\textbf{Proof}. 
Let $H:=\{x_1=0\}$ and let $X:=\{x_1=\cdots=x_k=0\}$.
Let $S = S^X \otimes_\K S_X$ with $S^X:=\mbox{Sym}^*(X^*)=\K[x_{k+1},\ldots,x_\ell]$ and 
$S_X:=\mbox{Sym}^*((V/X)^*)=\K[x_1,\ldots,x_k]$. Note that, by Proposition \ref{directsum}, 
$$
D(\A_X^e) \otimes_{\K} S^X =S \cdot D(\A_X^e),\ S \cdot D(\A_X^e) \oplus F = D(\A_X),
$$
where $F:=\oplus_{i={k+1}}^\ell S \partial_{x_i}$ is a free $S$-module. 
Since the localization is an exact functor and $D(\A)_p=D(\A_X)_p$ for $p=(x_1,\ldots,x_k)$ by Lemma \ref{localfunctor}, 
the surjectivity of $\rho^H$ 
shows that 
$(\rho_X^H)_p:D(\A_X)_p \rightarrow D(\A_X^H)_p$ is surjective. 
Let $C:=\coker \rho_X^H$, which satisfies $C_p=0$. It suffices to show that $C=0$.
Let $C^e$ be the essentialization of $C$, that is an $S_X$-module. Then $C=0 \iff C^e=0 \iff 
C^e_{p \cap S_X}=0$ since $C^e$ is an $S_X$-graded module. Thus we may replace $C^e$ by 
$C^e_{p \cap S_X}$. 
Since $\otimes_{(S_X)_{p \cap S_X}} S_p$ is a faithfully 
flat functor, and 
$0=C_p=C^e_{p \cap S_X}\otimes_{S_{p \cap S_X}} S_p$, it holds that $C^e=0$.  \owari
\medskip

\begin{lemma}
Let $H \in \A$. Assume that the Ziegler restriction $\pi^H:D_H(\A) \rightarrow 
D(\A^H,m^H)$ is surjective. Then the localization 
$\pi_X^H:D_H(\A_X) \rightarrow D(\A_X^H,m^H_X)$ is surjective for 
$X \in L(\A^H)$.
\label{surjlocpi}
\end{lemma}

\noindent
\textbf{Proof}. 
Let $H:=\{x_1=0\}$ and let $X:=\{x_1=\cdots=x_k=0\}$.
Let $S = S^X \otimes_\K S_X$ with $S^X:=\mbox{Sym}^*(X^*)$ and 
$S_X:=\mbox{Sym}^*((V/X)^*)$. Note that 
$$
D_H(\A_X^e) \otimes_{\K} S^X =S \cdot D_H(\A_X^e),\ S \cdot D_H(\A_X^e) \oplus F = D_H(\A_X),
$$
where $F:=\oplus_{i={k+1}}^\ell S \partial_{x_i}$ is a free $S$-module. 
For $p=(x_1,\ldots,x_k)$, Lemma \ref{localfunctor} and 
the definition of $D_H(\A)$ say that $D_H(\A_X)_p=D_H(\A)_p$. 
Since the localization is an exact functor, 
the surjectivity of $\pi^H$ shows that 
$(\pi_X^H)_p:D_H(\A_X)_p \rightarrow D(\A_X^H,m^H)_p$ is surjective. 
Let $C:=\coker \pi_X^H$, which satisfies $C_p=0$. It suffices to show that $C=0$.
Let $C^e$ be the essentialization of $C$, that is an $S_X$-module. Then $C=0 \iff C^e=0 \iff 
C^e_{p \cap S_X}=0$ since $C^e$ is an $S_X$-graded module. Thus we may replace $C^e$ by 
$C^e_{p \cap S_X}$. 
Since $\otimes_{(S_X)_{p \cap S_X}} S_p$ is a faithfully 
flat functor, and 
$0=C_p=C^e_{p \cap S_X}\otimes_{S_{p \cap S_X}} S_p$, it holds that $C^e=0$.  \owari
\medskip

For the projective dimension of $\A^H$, the following is fundamental, and we use it frequently without mentioning 
in the rest of this article.

\begin{lemma}
For $ H \in \A$, it holds that 
$$
\pd \A^H+1=\pd_S D(\A^H).
$$
\label{pdah}
\end{lemma}

\noindent
\textbf{Proof}.
Recall that $\pd \A^H:=\pd_{S^H} D(\A^H)$. By definition of depth, it holds that 
$$
\mbox{depth}_S D(\A^H)=\mbox{depth}_{S^H} D(\A^H)=:d.
$$
So by Auslander-Buchsbaum formula, it holds that 
$$
\pd \A^H=\ell-1-d=(\ell-d)-1=\pd_S D(\A^H)-1,
$$
which completes the proof.\owari
\medskip

\begin{lemma}
Let $(\A,m)$ be a multiarrangement such that 
$$
b_2(\A,m)=b_2(\A)+|\A|(|m|-|\A|).
$$
Then $\pd \A=\pd (\A,m)
$ unless $\pd (\A,m)=0$. In particular, 
$$
\pd \A \le1 \iff \pd (\A,m)\le 1.
$$ 
Moreover, the only case when $\pd \A\neq \pd (\A,m)$ is when 
$\pd \A=1>0=\pd (\A,m)$. Hence $\pd (\A,m) \le \pd \A$ 
under this assumption.
\label{pdeq}
\end{lemma}

\noindent
\textbf{Proof}. 
By the same argument as that of Theorem \ref{b2gen} with 
the given equality, there is a set of generators $\theta_1,\ldots,\theta_s$ together with $Q'\theta_E$ for $D(\A,m)$ such that $\theta_E, \theta_1,\ldots,\theta_s$ form a generator for $D(\A)$. 
Here $Q':=Q(\A,m)/Q(\A)$. Let $d_i:=\deg \theta_i$, $d:=\deg Q'+1$ and let 
$$
0 \rightarrow K  \stackrel{F}{\rightarrow} S[-1] \oplus (\oplus_{i=1}^s S[-d_i]) \rightarrow D(\A) \rightarrow 0
$$
and 
$$
0 \rightarrow K' \stackrel{F'}{\rightarrow} S[-d] \oplus (\oplus_{i=1}^s S[-d_i])\rightarrow D(\A,m) \rightarrow 0
$$
be free resolutions. It is clear that $K' 
\subset K$. Let us show that $K \subset K'$. Since $\theta_i \in D(\A,m)$ for 
$i=1,\ldots,s$, for every $0 =f\theta_E +\sum_{i=1}^s f_i \theta_i \in K$, it holds that 
$Q' \mid f$. Thus $K=K'$. So for $i \ge 2$, it holds that 
$$
\mbox{Ext}_S^i(D(\A),S) \simeq \mbox{Ext}_S^{i-1}(K,S) =\mbox{Ext}_S^{i-1}(K',S) \simeq \mbox{Ext}_S^i(D(\A,m),S).
$$
So $\pd \A \le 1 \iff \pd (\A,m) \le 1$, and they coincide if $\pd \A \ge 2$ or $\pd (\A,m) \ge 2$. 
The rest part follows from Theorem \ref{division} (4). \owari
\medskip

A trivial corollary of the above proof is the following.

\begin{cor}
Assume that 
$$
b_2(\A,m)=b_2(\A)+|\A|(|m|-|\A|).
$$
Then there is a commutative 
diagram of free resolutions of $D(\A)$ and 
$D(\A,m)$ of the following forms:
$$
\xymatrix{
0 \ar[r] &K  \ar[r]^-F &S[-1] \oplus (\oplus_{i=1}^s S[-d_i]) \ar[r]& D(\A) \ar[r] &0\\
0 \ar[r] &K \ar@{=}[u]\ar[r]^-{F'} &S[-d] \oplus (\oplus_{i=1}^s S[-d_i])\ar[r] \ar@<+7ex>@{^{(}->}[u]
\ar@<-5ex>@{=}[u]
&D(\A,m) \ar[r] \ar@{^{(}->}[u]&0
}
$$
Here $d:=|m|-|\A|+1$. 
In particular, for $i \ge 2$, it holds that 
$$
\mbox{Ext}_S^i(D(\A),S) \simeq \mbox{Ext}_S^i(D(\A,m),S).
$$
\label{kernel}
\end{cor}

To prove Theorem \ref{pisurj}, we need the following two lemmas.

\begin{lemma}
Assume that $\pi^H$ is surjective. Then the multi-$b_2$-equality holds for $(\A,H)$.
\label{lemma1}
\end{lemma}

\noindent
\textbf{Proof}. 
Note that every localization of a surjective $\pi^H$ is again surjective by Lemma \ref{surjlocpi}. 
So for every $X \in L_2(\A^H)$, the map $\pi_X^H:D_H(\A_X) 
\rightarrow D(\A_X^H,m^H_X)$ 
is surjective. Thus Theorem \ref{Y} shows that $\A$ is locally free in codimension three 
along $H$, which is equivalent to the multi-$b_2$-equality for $(\A,H)$ by Theorem \ref{Y}.\owari
\medskip

The following enables us to compare the second Betti numbers of $D(\A^H)$ 
with 
$D(\A^H,m^H)$. 

\begin{lemma}
The lower $b_2$-equality $b_2(\A^H,m^H)=b_2(\A^H) +(|\A^H|-1)(|\A|-|\A^H|-1)$ holds if and only if 
$D(\A^H)=D(\A^H,m^H)+S^H \rho^H(\theta_E)$.
\label{lemma2}
\end{lemma}

\noindent
\textbf{Proof}.
When $\ell \le 2$, then there is nothing to prove. 
Let $\ell=3$. Since $b_2(\A^H,m^H)=d_1d_2$ if $\exp(\A^H,m^H)
=(d_1,d_2)$, 
$\exp(\A^H)=(1,|\A^H|-1)$, Theorem \ref{adddel} shows that 
the lower $b_2$-equality is equivalent to 
$\exp(\A^H,m^H)=(|\A|-|\A^H|,|\A^H|-1)$. In this case, there is a basis $\theta_E,\theta$ for $D(\A^H)$ such that 
$f\theta_E,\theta$ also form 
a basis for $D(\A^H,m^H)$ for some polynomial $f$. Thus the statement is clear.

Assume that 
$\ell \ge 4$. 
First let us show the ``only if'' part. By Theorem \ref{b2gen}, 
there is a set of generators $\theta_2,\ldots,\theta_s$ for $D(\A^H,m^H)$ such that 
together with $\rho^H(\theta_E)$ form that for $D(\A^H)$, which completes the proof.

Next let us prove the ``if'' part. 
By Corollary 
\ref{kernel}, we have two free 
resolutions:
\begin{eqnarray*}
&0& \rightarrow K \rightarrow \oplus_{i=1}^s S^H[-d_i] \oplus S^H[-1] \rightarrow D(\A^H) \rightarrow 0,\\
&0& \rightarrow K \rightarrow \oplus_{i=1}^s S^H[-d_i] \oplus S^H[-|\A|+|\A^H|] \rightarrow D(\A^H,m^H) \rightarrow 0.
\end{eqnarray*}
Let $b_i^K$ and $b_i^F$ be the $i$-th Chern classes of the sheafified $K$ and $F$. Then 
by the functoriality of Chern classes and Theorem \ref{DS}, we have two equalities modulo $t^3$:
\begin{eqnarray*}
(1-|\A^H|t+b_2(\A)t^2)(1-b_1^Kt+b_2^Kt^2)&=&(1-b_1^Ft+b_2^Ft^2)(1-t),\\
(1-(|\A|-1)t+b_2(\A,m)t^2)(1-b_1^Kt+b_2^Kt^2)&=&(1-b_1^Ft+b_2^Ft^2)(1-(|\A|-|\A^H|)t).
\end{eqnarray*}
Thus 
\begin{eqnarray*}
b_1^F-b_1^K&=&|\A^H|-1,\\
b_2(\A)+|\A^H|b_1^K+b_2^K&=&b_2^F+b_1^F,\\
b_2(\A,m)+(|\A|-1)b_1^K+b_2^K&=&b_2^F+b_1^F(|\A|-|\A^H|).
\end{eqnarray*}
From these equations, we have
$$
b_2(\A^H,m^H)-b_2(\A^H)=(b_1^F-b_1^K)(|\A|-|\A^H|-1)=
(|\A^H|-1)(|\A|-|\A^H|-1),
$$
which is nothing but the lower $b_2$-equality.\owari
\medskip

Now we can show one direction of Theorem \ref{pisurj} without the assumption NMPD. 

\begin{theorem}
Assume that $\rho^H$ and $\pi^H$ are both surjective for some $H \in \A$. Then 
the $b_2$-equality holds for $(\A,H)$.
\label{surjkarab2}
\end{theorem}

\noindent
\textbf{Proof}. 
Note that, by Lemma \ref{lemma1}, the surjectivity of 
$\pi^H$ implies the upper $b_2$-equality. Thus 
by Proposition \ref{b2prop} (3), it suffices to show the lower $b_2$-equality holds. 

Recall that $D(\A)=S\theta_E \oplus D_H(\A)$. Thus by the definitions and surjectivity of $\rho^H$ and $\pi^H=\rho^H|_{D_H(\A)}$, it holds that 
$$D(\A^H)=\mbox{Im}(\rho^H)=\mbox{Im}(\pi^H) + S^H \rho^H(\theta_E)=
D(\A^H,m^H)+S^H \rho^H(\theta_E)$$ by the assumptions and definitions. Thus Lemma 
\ref{lemma2} shows that the lower $b_2$-equality holds. \owari
\medskip

On the implication from the (multi-)$b_2$-equality to surjectivity, the following proposition 
on the cohomology vanishing is the key. 

\begin{prop}
Let $E$ be a coherent sheaf on $\P^{n}$ with $n \ge k+2$. Assume that $E\simeq 
\widetilde{M}$ for some $S:=H_*^0(\mathcal{O})$-module $M$, and 
$\pd_S M =j \le k$. Then 
$$
H^i_*(E):=\oplus_{d \in \Z} H^i(E(d),\P^n)=0
$$
for $1 \le i \le k+1-j$. 
\label{pd0}
\end{prop}

\noindent
\textbf{Proof}. 
Induction on $j$. If $j=0$, then the statement is clear since all line bundles have 
zero middle cohomologies. Assume that $j>0$. Let 
$$
0 \rightarrow K \rightarrow F_0 \rightarrow E \rightarrow 0
$$
be the last part of the free resolution of $E$ of length $j$. Since $\pd E =j$, 
it holds that $\pd K=j-1$. By induction hypothesis, 
$H^i_*(K)=0
$ for $1 \le i \le k+2-j<n$. Thus the long exact sequence of the above short exact sequence 
shows that 
$H^i_*(E)=0
$ for $1 \le i \le k+1-j$.\owari
\medskip

Now we completed the preparation for the proof of three main results on 
surjectivity. 
\medskip

\noindent
\textbf{Proof of Theorem \ref{pisurj}}. 
First let us show the ``if'' part, i.e., let us assume that 
$\pi^H$ is surjective. Then the statement follows directly from Lemma \ref{lemma1}. 
Next let us show the ``only if'' part. Assume that $b_2^0(\A)=b_2(\A^H,m^H)$ and 
$\A$ is NMPD along $H$.
Let us show by induction on $\ell$. When $\ell \le 3$, the $b_2$-equality implies that 
$\pi^H$ is surjective by Theorem \ref{Y}. Assume that $\ell \ge 4$. 
Consider $\A_X$ for $0 \neq X \in L_s(\A^H)$ with $s \ge 2$. First let us show that 
$\pi_X^H:D(\A_X^e) \rightarrow D((\A_X^H)^e)$ is surjective. 
Since $\A$ is NMPD along $H$, $\A_X^e$ is also 
NMPD along $H$. Also, the $b_2$-equality holds for $(\A_X^e,H)$ by 
Proposition \ref{b2local}. 
Thus induction hypothesis confirms that $\pi_X^H$ is surjective, and so is its localization. Hence we have 
the exact sequence 
$$
0 \rightarrow \widetilde{D_H(\A)}
\stackrel{\cdot \alpha_H}{\rightarrow} \widetilde{D_H(\A)}
\stackrel{\pi}{\rightarrow} \widetilde{D(\A^H,m^H)} \rightarrow 0.
$$
Thus by Proposition \ref{pd0} and the fact that 
$$
D(\A)=S[-1] \oplus D_H(\A),$$ 
we can apply 
Proposition 
\ref{globalsection} to the cohomology long 
exact sequence  
$$
0 \rightarrow D_H(\A) \stackrel{\cdot \alpha_H}{\rightarrow} D_H(\A) \stackrel{\pi}{\rightarrow} D(\A^H,m^H) \rightarrow H^1_*(\widetilde{D_H(\A)}).
$$
By Proposition \ref{pd0} and NMPD, it holds that $H^1_*(\widetilde{D_H(\A)})=0$, 
thus $\pi^H$ is surjective. 
\owari
\medskip

\noindent
\textbf{Proof of Theorem \ref{surj}}. 
(1)\,\,
The ``only if'' part follows from Lemma \ref{surjloc}. 
Let us show the ``if'' part. The statement is true if 
$\ell \le 3$  by 
the 
assumption on $\rho^H$. Assume that the statement is true up to $\ell-1 \ge 3$. 
Consider $\B:=(\A_X')^e$ for $0 \neq X \in L_s(\A^H)$ with $s \ge 2$. Let us show that 
$\rho_X^H:D(\A_X^e) \rightarrow D((\A_X^H)^e)$ is surjective. 
Since $\A'$ is NMPD along $H$, $\B$ is also 
NMPD along $H$. Also, $\rho_X^H$ is
surjective if $s=2$ by the assumption. 
Thus induction hypothesis confirms that $\rho_X^H$ is surjective, and so is its localization. Hence we have 
the exact sequence 
$$
0 \rightarrow \widetilde{D(\A')}
\stackrel{\cdot \alpha_H}{\rightarrow} \widetilde{D(\A)}
\stackrel{\rho^H}{\rightarrow} \widetilde{D(\A^H)} \rightarrow 0.
$$
Taking the global section combined with Proposition \ref{globalsection}, we have 
$$
0 \rightarrow D(\A') \stackrel{\cdot \alpha_H}{\rightarrow} D(\A) \stackrel{\rho^H}{\rightarrow} D(\A^H) \rightarrow H^1_*(\widetilde{D(\A')}).
$$
Since $\pd \A'$ is not maximal, Proposition \ref{pd0} shows that 
$H^1_*(\widetilde{D(\A')})=0$. 
Thus $\rho^H$ is surjective.

(2)\,\,The ``if'' part follows from Theorem \ref{surjkarab2}. 
Next let us show the ``only if'' part. 
Since the $b_2$-equality holds, the multi-$b_2$-equality holds too. 
Hence $\pi^H$ is 
surjective by Theorem \ref{pisurj}. Now 
Proposition \ref{pirho} completes the proof.

(3)\,\, 
Assume that 
$\A'$ is NMPD along $H$. 
By Theorems \ref{Y} and \ref{division}, the $b_2$-equality implies that 
$\A$ is locally surjective in codimension three. Thus (1) shows that 
$\rho^H$ is surjective. \owari
%
\medskip

\begin{example}
Let $\A$ be defined as
$$
x_1x_2x_3x_4(x_1+x_2+x_3+x_4)=0
$$
and let $H:x_1=0$. Then $\A'$ is free, the $b_2$-equality holds for $(\A,H)$ and 
$\pd \A$ is maximal. Since $D_H(\A)$ has no 
degree one 
derivations and $\rho^H(\theta_E) \in 
D(\A^H,m^H)=D(\A^H)$, it holds that $\pi^H$ is not surjective. Thus 
the assumption that $\pd \A$ is not maximal is necessary in Theorem \ref{surj}.
\label{surjce}
\end{example}

\noindent
\textbf{Proof of Theorem \ref{freesurj}}. 
Since $\pd\A'=0$, it holds that 
$$
H^1_*(\widetilde{D(\A')})=0
$$
when $\ell \ge 3$. 
We prove by induction on $\ell \ge 1$. When $\ell \le 2$ there is nothing to 
show. Assume that $\ell \ge 3$. Since $\A'_X$ is free too by 
Proposition \ref{pdlow}, the induction hypothesis gives 
$$
0 \rightarrow \widetilde{D(\A')}
\stackrel{\cdot \alpha_H}{\rightarrow} \widetilde{D(\A)}
\stackrel{\rho^H}{\rightarrow} \widetilde{D(\A^H)} \rightarrow 0.
$$
Taking the global section combined with Proposition \ref{globalsection}, we have 
$$
0 \rightarrow D(\A') \stackrel{\cdot \alpha_H}{\rightarrow} D(\A) \stackrel{\rho^H}{\rightarrow} D(\A^H) \rightarrow H^1_*(\widetilde{D(\A')})=0,
$$
thus $\rho^H$ is surjective. 
\owari
\medskip


Now we can show the addition-deletion theorems.
\medskip

\noindent
\textbf{Proof of Theorem \ref{pdadd}}.
%
%


(1)\,\,
By Theorem \ref{surj} (1), we have the exact sequence 
$$
0 \rightarrow D(\A') \stackrel{\cdot \alpha_H}{\rightarrow} D(\A) \stackrel{\rho^H} \rightarrow 
D(\A^H) \rightarrow 0.
$$
Since $\pd \A^H=k
=\pd \A'$, we have 
\begin{eqnarray*}
\mbox{Ext}^i(D(\A'),S)&=&0\ (i >k),\\
\mbox{Ext}^i(D(\A^H),S)&=&0\ (i >k+1).
\end{eqnarray*}
Hence 
$$
\mbox{Ext}^i(D(\A),S)=0\ (i >k+1).
$$
Thus $\pd \A\le k+1$.

(2)\,\,
Let $\pd \A^H=s>k$. Thus 
\begin{eqnarray*}
\mbox{Ext}^k(D(\A'),S)&\neq&0,\\
\mbox{Ext}^i(D(\A'),S)&=&0\ (i >k),\\
\mbox{Extt}^i(D(\A^H),S)&=&0\ (i >s+1),\\
\mbox{Ext}^{s+1}(D(\A^H),S)&\neq&0.
\end{eqnarray*}
Hence 
$$
\mbox{Ext}^i(D(\A),S)=0\ (i \ge s+2)
$$
and 
$$
0 \neq \mbox{Ext}^{s+1} (D(\A^H),S) \simeq 
\mbox{Ext}^{s+1} (D(\A),S).
$$ 
Hence $\pd \A=s+1=\pd \A^H+1$.

(3)\,\,
Let $\pd \A^H=s<k$. Then 
\begin{eqnarray*}
\mbox{Ext}^k(D(\A'),S)&\neq&0,\\
\mbox{Ext}^i(D(\A'),S)&=&0\ (i >k),\\
\mbox{Extt}^i(D(\A^H),S)&=&0\ (i >s+1 \le k).
\end{eqnarray*}
Hence 
$$
\mbox{Ext}^i(D(\A),S)=0\ (i \ge k+1)
$$
and 
$
0 \neq \mbox{Ext}^{k} (D(\A),S)$ since 
$\mbox{Ext}^{k} (D(\A),S) \rightarrow  \mbox{Ext}^{k} (D(\A'),S) \rightarrow 0$ and 
$\mbox{Ext}^{k} (D(\A'),S) \neq 0$ since $\pd_S \A'=k$.
Hence $\pd \A=k=\pd \A'>\pd \A^H$.
\owari
\medskip

\noindent
\textbf{Proof of Theorem \ref{pddel}}.

By Theorem \ref{surj} (1), we have the exact sequence 
$$
0 \rightarrow D(\A') \stackrel{\cdot \alpha_H}{\rightarrow} D(\A) \stackrel{\rho^H} \rightarrow 
D(\A^H) \rightarrow 0.
$$

(1)\,\,
Clear since $\mbox{Ext}^{i \ge k+1}(D(\A),S)=\mbox{Ext}^i(D(\A^H),S)=0$. 

(2)\,\,
Let $\pd \A^H=s>k-1$. Thus 
\begin{eqnarray*}
\mbox{Ext}^k(D(\A),S)&\neq&0,\\
\mbox{Ext}^i(D(\A),S)&=&0\ (i >k),\\
\mbox{Ext}^i(D(\A^H),S)&=&0\ (i >s+1),\\
\mbox{Ext}^{s+1}(D(\A^H),S)&\neq&0.
\end{eqnarray*}
Hence 
$$
\mbox{Ext}^i(D(\A'),S)=0\ (i \ge s+1)
$$
and 
$$
0 \neq \mbox{Ext}^{s+1} (D(\A^H),S) \simeq 
\mbox{Ext}^{s} (D(\A'),S).
$$ 
Hence $\pd \A'=s=\pd \A^H$.

(3)\,\,
Let $\pd \A^H=s<k-1$. Thus 
\begin{eqnarray*}
\mbox{Ext}^k(D(\A),S)&\neq&0,\\
\mbox{Ext}^i(D(\A),S)&=&0\ (i >k),\\
\mbox{Ext}^i(D(\A^H),S)&=&0\ (i >s+1\le k-1).
\end{eqnarray*}
Hence 
$$
\mbox{Ext}^i(D(\A'),S)=0\ (i \ge k+1)
$$
and 
$$
0 \neq \mbox{Ext}^{k} (D(\A'),S) \simeq 
\mbox{Ext}^{k} (D(\A),S) \neq 0.
$$ 
Hence $\pd \A'=\pd \A>\pd \A^H+1$.\owari

\medskip

\noindent
\textbf{Proof of Theorem \ref{pdrest}}. 
(1)\,\,
By Theorem \ref{surj} (1), we have the exact sequence 
$$
0 \rightarrow D(\A') \stackrel{\cdot \alpha_H}{\rightarrow} D(\A) \stackrel{\rho^H} \rightarrow 
D(\A^H) \rightarrow 0.
$$
Let $\pd \A=s<k$. Thus 
\begin{eqnarray*}
\mbox{Ext}_S^k(D(\A'),S)&\neq&0,\\
\mbox{Ext}_S^i(D(\A'),S)&=&0\ (i >k),\\
\mbox{Ext}_S^i(D(\A),S)&=&0\ (i >s),\\
\mbox{Ext}_S^{s}(D(\A),S)&\neq&0.
\end{eqnarray*}
Hence 
$$
\mbox{Ext}_S^i(D(\A^H),S)=0\ (i >k+1)
$$
and 
$$
\mbox{Ext}^{k+1}_S(D(\A^H),S) \simeq 
\mbox{Ext}^k_S(D(\A'),S)\neq 0.
$$ 
Thus $\pd \A^H=k$. 

(2)\,\,
Since $\pd \A=k
=\pd A'$, we have 
\begin{eqnarray*}
\mbox{Ext}^i_S(D(\A'),S)&=&0\ (i >k),\\
\mbox{Ext}^i_S(D(\A),S)&=&0\ (i >k).
\end{eqnarray*}
Hence 
$$
\mbox{Ext}^i(D(\A^H),S)=0\ (i >k+1).
$$
Thus $\pd \A^H\le k$.

(3)\,\,
Let $\pd \A=s>k$. Thus 
\begin{eqnarray*}
\mbox{Ext}^k_S(D(\A'),S)&\neq&0,\\
\mbox{Ext}^i_S(D(\A'),S)&=&0\ (i >k),\\
\mbox{Ext}^i_S(D(\A),S)&=&0\ (i >s),\\
\mbox{Ext}^s_S(D(\A),S)&\neq&0.
\end{eqnarray*}
Hence 
$$
\mbox{Ext}^i_S(D(\A^H),S)=0\ (i \ge s+1).
$$
Since $\mbox{Ext}^{s} (D(\A^H),S)=0$ implies that 
$\mbox{Ext}^{s} (D(\A),S)=0$, we have 
$\mbox{Ext}^{s} (D(\A^H),S)\neq 0$. 
Hence $\pd \A^H=s-1=\pd \A-1$. \owari
%
\medskip

%

Let us see examples to which we can apply Theorems \ref{pdadd}, 
\ref{pddel} and 
\ref{pdrest}.

\begin{example}
(1)\,\,
Let $$
Q(\A)=xyz(x+y+z)
$$
and $\A \ni H:=\{x+y+z=0\}$. 
Then $\A'$ is free, thus it is NMPD. Thus Theorem \ref{freesurj} shows that $\rho^H$ is surjective in codimension three. Also, $\A^H$ is free since it is in $\K^2$. Thus we may apply Theorem \ref{pdadd} (1) to know that 
$$
\pd \A'=\pd \A^H=0,\ 1=\pd \A\le 0+1=1.
$$
This can be regarded to the example corresponding to Theorem \ref{pddel} (2). 

(2)\,\,
Let $$
Q(\A)=xyzw(x+y+z+w)
$$
and $\A \ni H:=\{x+y+z+w=0\}$. 
Then $\pd \A'=0$ is free, thus it is NMPD. Thus Theorem \ref{freesurj} shows that $\rho^H$ is surjective in codimension three. Also, $\pd \A^H=1>\pd \A'=0$ by (1). Thus we may apply Theorem \ref{pdadd} (2) to know that 
$$
\pd \A'=0,\ \pd \A^H=1,\ \pd \A =\pd \A^H+1=2.
$$
This can be regarded to the example corresponding to Theorem \ref{pdrest} (3). 

(3)\,\,
Let $$
Q(\A)=\prod_{i=1}^5 x_i \prod (x_1 \pm x_2 \pm x_3 \pm  x_4 \pm x_5).
$$
in $\R^5$ ,
and $\A \ni H:=\{x_1-x_2-x_3-x_4-x_5=0\}$. Edelman and Reiner showed in \cite{ER} that 
$\A$ is free, but $\A^H$ and $\A'$ are both not free. By Theorem \ref{SPOG}, we know that $\pd \A'=1$, thus it is NMPD by Theorem \ref{NMPDcond}. Also, since $\A$ is NMPD, and the $b_2$-equality holds for $(\A,H)$, Theorem \ref{surj} (2) shows that $\rho^H$ is surjective in codimension three. Hence we may apply Theorem \ref{pdrest} (1) to know that $\pd \A^H=\pd \A'=1$. 


\end{example}

\section{Proof of Theorem \ref{PDCcombin}}

To prove Theorem \ref{PDCcombin}, first let 
us prove the following:

\begin{prop}
Assume that $\pd \A'=\ell-2$, i.e., maximal. 
If the $b_2$-equality holds for $(\A,H)$ and $\pd \A_X <\codim X-2$ for all $X \in L(\A) \setminus \{0\}$, then $\pd \A=\ell-2$. 
\label{maxb2}
\end{prop}

\noindent
\textbf{Proof}. 
Assume that $\pd \A< 
\ell-2$. Then by Theorem \ref{surj} (2), 
both $\pi^H$ and 
$\rho^H$ are surjective. Note that $\mbox{Ext}^{\ell-1}_{S}(D(\A^H),S)=0$. 
Thus by the exact sequence 
$$
0 \rightarrow D(\A') \stackrel{
\cdot \alpha_H}{\rightarrow} D(\A) \stackrel{\rho^H}{\rightarrow} D(\A^H) \rightarrow 0
$$
and its long exact sequence of Ext's, we have $\mbox{Ext}^{\ell-2}(D(\A'),S)=0$, a contradiction.\owari
\medskip

\begin{example}
We say that $\A$ is \textbf{generic} if $\A$ is irreducible and 
$\codim X=|\A_X|$ for all $X \in L(\A) \setminus \{0\}$. It is known that 
$\pd \A=\ell-2$, see Corollary 4.4.3 in \cite{RT} for example. We prove this by using Theorems \ref{pdadd} and 
Proposition \ref{maxb2}. The statement is trivial if $\ell \le 3$. Assume that the statement is true up to $\ell-1 \ge 3$. 
We also prove by induction on $|\A|$. The smallest case is $|\A|=\ell+1$, that is of the form 
$$
(x_1+\cdots+x_\ell)\prod_{i=1}^\ell x_i =0.
$$
Let $H:x_1+\cdots+x_\ell=0$. 
It is clear that $\A'=\A \setminus \{H\}$ is free. 
Since 
$$
\A^H:(x_1+\cdots+x_{\ell-1})\prod_{i=1}^{\ell-1} x_i =0,
$$
which is generic, 
the induction hypothesis implies that $\pd \A^H=\ell-3>0=\pd \A'$. 
Thus Theorem \ref{pdadd} (2) shows that $\pd \A=\ell-2$. 

Now let us assume that $\A'$ is generic and $\pd \A'=\ell-2$. 
Let $\A:=\A' \cup \{H\}$ is generic too. Let us prove that $\pd \A=\ell-2$. 
Since $\A$ is locally free and the $b_2$-equality holds by genericity of $\A$ and Proposition \ref{b2local}, 
Proposition 
\ref{maxb2} shows that $\pd \A=\ell-2$. 
\end{example}

To prove Theorem \ref{PDCcombin}, let us complete the proof 
of Theorem \ref{add1}.
\medskip

\noindent
\textbf{Proof of Theorem \ref{add1}}. 
By Theorem \ref{freesurj}, in this case we have the Euler exact sequence. Thus Ext-long exact sequence 
shows the first part. 
(1) and (2) follow immediately from Theorems \ref{adddel} and \ref{division}. \owari
\medskip

Now we can prove Theorem \ref{PDCcombin}.
\medskip

\noindent
\textbf{Proof of Theorem \ref{PDCcombin}}.
First, if $\chi_0(\A;t)$ is irreducible over $\Z$, then Terao's factorization shows that $\A$ is not free. 
Thus $\pd \A=1$ if $\ell=3$. Also, $\chi_0(\A;|\A^H|-1) \neq 0$ 
implies that $\A$ is not free if $\A \setminus \{
H\}$ is free by Theorems \ref{adddel} and \ref{division}. Thus $\pd \A=1$ for such $\A$ with 
$\ell
=3$, which are combinatorial. 
Since the all operations of $\PDC$ 
including the above are combinatorial, it is clear by Theorem \ref{pdadd}, \ref{add1}, 
Propositions \ref{maxb2} and \ref{CScond}. \owari
\medskip

When we apply Theorem \ref{PDCcombin}, we have to check whether $(\A,H) \in \mbox{CS}_3$ or not. So 
let us introduce some conditions for 
$\A \in \mbox{CS}_3$. 

\begin{prop}
Let $\A$ be an arrangement in $\K^3$ and 
let $H \in \A,\ \A':=\A \setminus \{H\}$. Then 
$(\A,H) \in \mbox{CS}_3$ if 
there exists at most 
one $X_0 \in \A^H$ such that $m^H(X_0) >1$.
\label{CSsuff}
\end{prop}

\noindent
\textbf{Proof}.
Assume that 
there exists at most 
one $X_0 \in \A^H$ such that $m^H(X_0) =k+1>1$.
Since the proof is the same we may assume that there exist such $X_0$, defined by $y=z=0$ and $H:z=0$. 
We may also assume that 
$\{x=0\} 
\in \A$. Then 
$$
Q(\A^H)=xf(x,y),\ 
Q(\A^H,m^H)=x^{k+1}f(x,y)
$$
with $f(0,y) \neq 0$. Then 
$$
D(\A^H)=\langle 
\rho^H(\theta_E), f\partial_y
\rangle
$$
by Theorem \ref{saito}. 
Then we may express
$$
Q(\A)=zQ_1Q_2,
$$
where 
$Q_1|_{z=0}=x^{k+1}$ and 
$Q_2=Q_2(x,y,z)$ with $Q_2(0,y,0) \neq 0$ and 
$Q_2(x,y,0)=f(x,y)$. Thus 
$$
Q_2\partial_{y} 
\in D(\A)$$
and $\rho^H(Q_2 \partial_y)=f\partial_y$, which shows that 
$\rho^H$ is surjective.\owari
\medskip

Summarizing, we have the following:

\begin{prop}
Let $\A$ be an arrangement in $\K^3$ and 
let $H \in \A,\ \A':=\A \setminus \{H\}$. Then 
$(\A,H) \in \mbox{CS}_3$ if one of the following conditions holds true:
\begin{itemize}
\item[(1)]
$\A' \in \mathcal{DF}$. Namely, there is $L \in \A'$ such that the $b_2$-equality holds for $(\A',L)$. 
\item[(2)]
There exists at most 
one $X_0 \in \A^H$ such that $m^H(X_0) >1$.
\item[(3)]
There is the $b_2$-equality for $(\A,H)$.
\end{itemize}\label{CScond}
\end{prop}

\noindent
\textbf{Proof}. Clear by Theorems \ref{SFDF}, \ref{freesurj}, 
\ref{Y} and Proposition \ref{CSsuff}. \owari
\medskip








\section{Applications, variants and examples}

Let us collect several applications, variants and examples 
of the results in the previous sections.

\subsection{Yoshinaga-type result for projective dimensions}

We can formulate the projective dimensional version of Yoshinaga's criterion (Theorem \ref{Y}) 
as follows:

\begin{theorem}
Let $H \in \A$. Assume that 
$\A$ is locally free along $H$, i.e., $\A_X$ is free for all $X \in L(\A^H) \setminus \{0\}$, and 
$\pd (\A^H,m^H)$ is not maximal. 
\begin{itemize}
\item[(1)]
Then $\pi^H$ is surjective and $\pd \A \le \pd (\A^H,m^H)$. 
\item[(2)]
Assume futher that the lower $b_2$-equality holds for $(\A,H)$ and 
$\pd \A^H$ is not maximal. Then 
$\pd \A \le \pd \A^H$.
\end{itemize}
\label{Ypd}
\end{theorem}

To prove Theorem \ref{Ypd} let us recall the following from \cite{Y1}.

\begin{prop}[\cite{Y1}, Theorem 2.3]
Let $E$ be 
a reflexive sheaf on $\P^n$ such that there are at most finitely many points $p \in \P^n$ such that $E_p$ is not free. Then 
$$
H^1(E(d))=0
$$
for all $d<<0$.
\label{YH10}
\end{prop}

\noindent
\textbf{Proof of Theorem \ref{Ypd}}. 
Let $\widetilde{D_H(\A)}$ be a reflexive sheaf on $\P^{\ell-1}$. Since 
$\A$ is locally free along $H$, we have two exact sequences
\begin{eqnarray*}
0 &\rightarrow& \widetilde{D_H(\A)} \stackrel{\cdot \alpha_H}
{\rightarrow} \widetilde{D_H(\A)} \stackrel{\pi^H}{\rightarrow} 
\widetilde{D(\A^H,m^H)} \rightarrow 0,\\
0 &\rightarrow& \widetilde{D_H(\A)} \stackrel{\cdot \alpha_H}
{\rightarrow} \widetilde{D_H(\A)} \stackrel{\pi^H}{\rightarrow} 
\widetilde{D_H(\A)}|_H \rightarrow 0,
\end{eqnarray*}
and it holds that 
$$
\widetilde{D(\A^H,m^H)} \simeq \widetilde{D_H(\A)}|_H.
$$
Since $\pd (\A^H,m^H)$ is not maximal, 
Proposition \ref{pd0} shows that 
$H^1_*(\widetilde{D_H(\A)}|_H)=0$. So the cohomology long exact sequence shows that 
$$
H^1(\widetilde{D_H(\A)}(d-1))
\stackrel{\cdot 
\alpha_H}{\rightarrow} H^1(\widetilde{D_H(\A)}(d))
$$
is surjective for all $d$. Therefore, 
Proposition \ref{YH10} confirms that 
$H^1_*(\widetilde{D_H(\A))})=0$. 
Hence we have a exact sequence 
$$
0 \rightarrow D_H(\A) 
\stackrel{\cdot 
\alpha_H}{\rightarrow} D_H(\A)
\stackrel{\pi^H}{\rightarrow} D(\A^H,m^H)
\rightarrow 0,
$$
which shows that 
$\pi^H$ is surjective.
Let $\pd \A^H=k$. Taking the Ext-long exact sequence, we have the surjection 
$\cdot \alpha_H:\mbox{Ext}^{k+1}(D_H(\A),S) \rightarrow \mbox{Ext}^{k+1}(D_H(\A),S)$ and 
$\mbox{Ext}^{i \ge k+2}(D_H(\A),S)=0$. 
Since $\cdot  \alpha_H$ is a surjection between the same $S$-graded module, $\alpha_H \cdot 
\mbox{Ext}_S^{k+1}(D(\A),S)=\mbox{Ext}_S^{k+1}(D(\A),S)$ shows that $\mbox{Ext}_S^{k+1}(D(\A),S)=0$, i.e., 
$\pd \A \le \pd (\A^H,m^H)$.
(2) follows immediately from (1) and Lemma \ref{pdeq}. 
\owari 
\medskip

\begin{example}
Let $\A$ be the Edelman-Reiner's example in \cite{ER} defined by 
$$
\prod_{i=1}^5 x_i\prod(x_1\pm x_2\pm x_3\pm x_4\pm x_5)=0
$$
in $\R^5$. Let 
$H:x_1-x_2-x_3-x_4-x_5=0$. 
Since 
$$
\chi(\A^H;t)=(t-1)(t-4)(t^2-10t+26), 
$$
Theorem \ref{Teraofactorization} shows that $\A^H$ is not free. 

First, note that 
$$
b_2(\A)=170,\ 
b_2(\A^H)=80,\ |\A|=21,\ |\A^H|=15.
$$
Thus $$
170=80+15(21-15),
$$
and the $b_2$-equality holds. Also, we can check that $\A$ is locally free along $H$ 
by a direct 
computation. Thus we can apply Theorem \ref{Ypd}. To do that, let us check $\pd \A^H$. 
Since we can show that $\A^H \cup\{x_1-x_2=0\}$ is free with exponents $(1,5,5,5)$, 
Theorem \ref{SPOG} shows that $\pd \A^H=1$, and Theorem \ref{Ypd} shows that 
$\pd \A \le \pd \A^H=1$. In fact, $\pd \A=0$, see 
Example \ref{ERex} for details. 
\label{ERex1}
\end{example}

\subsection{Division theorem for projective dimensions}

From a new viewpoint shown in this article, 
a generalization of Theorem \ref{division} 
to all projective dimensional cases can be given as follows:

\begin{theorem}[Division theorem for projective dimensions]
Assume that the $b_2$-equality holds for $(\A,H)$, and $\A$ is NMPD along $H$. 
Then 
$\ell-2>\pd \A'=\pd \A^H \ge \pd \A$.
\label{main}
\end{theorem}

\noindent
\textbf{Proof}. 
If $\pd \A^H=0$, then this is nothing but
Theorem \ref{division} asserting that 
$\pd \A=\pd \A^H=0$. Thus we may assume that $\pd \A^H>0$. 
First assume that $\A$ is NMPD and $\pd \A'$ is not maximal. 
Then by Theorem \ref{surj} (2) and 
Proposition \ref{pirho},
both $\pi^H$ and $\rho^H$ are 
surjective. 
First let us 
use the Ziegler 
restriction:
$$
0 \rightarrow D_H(\A) \stackrel{\cdot \alpha_H}{\rightarrow} D_H(\A) \stackrel{\pi^H} \rightarrow 
D(\A^H,m^H) \rightarrow 0
$$
which is right exact as shown above. 
Since $\pd \A^H=k \ge 1$ as assumed, 
Lemma \ref{pdeq} implies that $\pd (\A^H,m^H) \le k$. 
Thus $\mbox{Ext}_S^i(D(\A^H,m^H),S)=0$ 
for $i \ge k+2$. 
Hence we have surjections 
$$\cdot \alpha_H:\mbox{Ext}^{k+1}(D_H(\A),S)
\rightarrow \mbox{Ext}^{k+1}(D_H(\A),S).
$$
Since this is a surjection between the same $S$-graded module, $\alpha_H \cdot 
\mbox{Ext}_S^{k+1}(D(\A),S)=\mbox{Ext}_S^{k+1}(D(\A),S)$ shows that $\mbox{Ext}_S^{k+1}(D(\A),S)=0$. Hence 
$\mbox{Ext}_S^{k+1}(D_H(\A),S)=0$ because 
$$
D(\A)=S\theta_E \oplus D_H(\A).
$$
Applying the same argument to $i>k+1$, we have 
$\pd \A \le k=\pd \A^H$. 
Again applying the Ext-functor to the Euler 
restriction
$$
0 \rightarrow D(\A') \stackrel{\cdot \alpha_H}{\rightarrow} D(\A) \stackrel{\rho^H} \rightarrow 
D(\A^H) \rightarrow 0
$$
which is right exact as shown above,
it is clear that $\mbox{Ext}^{k}_S(D(\A'),S)\neq 0$ since 
$\mbox{Ext}^{k+1}_S(D(\A^H),S) \neq 0$. 
Because $\mbox{Ext}^{i \ge k+1}_S(D(\A'),S)=0$, 
it holds that $\pd \A'=\pd \A^H=k\ge \pd \A$. 

Next assume that $\pd \A'=\ell-2$ and $\A$ 
is NMPD. Then 
Theorem \ref{surj} (2) shows that both $\pi^H$ and $\rho^H$ are surjective. 
Hence the same argument as the above shows that $\mbox{Ext}^{\ell-2}_S(D(\A'),S)=0$, 
a contradiction. \owari
%
%
\medskip

Theorem \ref{main} gives us a rough estimate of 
projective dimensions when the $b_2$-equality holds, and also show that 
the $b_2$-equality holds rarely. The case when $\pd \A^H=0$ in Theorem \ref{main} is the division theorem (Theorem \ref{division}).

\subsection{Condition for NMPD}

In our main results, the key is the assumption ``NMPD''. As we have seen, it is not easy to check whether $\A$ is 
NMPD along $H$. However, the following gives us several ways to check MNPD.

\begin{theorem}
$\A$ is NMPD along $H$ if one of the following is satisfied:
\begin{itemize}
\item[(1)]
$\pd \A \le 1$.
\item[(2)]
$\pd \A_X \le 1$ for all $X \in L(\A^H) \setminus \{0\}$ and 
$\pd \A <\ell-2$. In particular, 
$\A$ is NMPD along $H$ if $\A$ is locally free and $\pd \A <\ell-2$.
\item[(3)]
$1<\pd \A=k<\ell-2$ and 
$\pd \A_X < s-2$ for all $X \in L_{s-1}(\A^H)$ with $
2\le s-2 \le k$.
\end{itemize}
\label{NMPDcond}
\end{theorem}

\noindent
\textbf{Proof}. (1) and (2) are clear 
from the definitions. For (3), apply Lemma \ref{pdlow}, Theorems \ref{main} and 
\ref{division}.\owari
\medskip

\subsection{Sheaf exact sequences}

There are some conditions for $\A$ to be NMPD shown Theorem \ref{NMPDcond}.
For the sheaves, we have the following.

\begin{cor}
Assume that the $b_2$-equality holds for $(\A,H)$, then 
\begin{itemize}
\item[(1)]
$$
0 \rightarrow \widetilde{D(\A')} 
\stackrel{\cdot \alpha_H}{\rightarrow}
\widetilde{D(\A)} 
\stackrel{\widetilde{\rho^H}}{\rightarrow}
\widetilde{D(\A^H)} 
\rightarrow 0
$$
is exact if $\A'$ is NMPD, and 
\item[(2)]
$$
0 \rightarrow \widetilde{D_H(\A)} 
\stackrel{\cdot \alpha_H}{\rightarrow}
\widetilde{D_H(\A)} 
\stackrel{\widetilde{\pi^H}}{\rightarrow}
\widetilde{D(\A^H,m^H)} 
\rightarrow 0.
$$
is exact if $\A$ is NMPD.
\end{itemize}
\label{sheaves}
\end{cor}

\noindent
\textbf{Proof}.
Note that $\pd \A_X$ is not maximal for any $X \in L(\A) \setminus \{0\}$. Thus 
the proofs of Theorems \ref{pdadd} and \ref{pddel} show that 
both $\rho_X^H$ and $\pi^H_X$ are 
surjective for $X \in L(\A^H) \setminus \{0\}$, which 
completes the proof.\owari
\medskip

\subsection{Surjectivity and freeness}
In this article the surjectivity of $\rho^H$ and $\pi^H$ played key roles. 
We can investigate them for details in this subsection. 
First, free surjection theorem and our new addition-deletion theorem give 
a necessary condition for $\A'$ to be free. 

\begin{cor}
Let $H \in \A$. Then $\A':=\A \setminus 
\{H\}$ is free only if $\rho^H$ is surjective and either (a) $\pd \A=\pd \A^H=0$, or (b) 
$\pd \A=\pd \A^H+1$.
\label{delfree}
\end{cor}

\noindent
\textbf{Proof}. 
Immediate from Theorems \ref{freesurj} and \ref{pddel}.
\owari
\medskip


When $\A$ is free, whether $\rho^H$ is 
surjective or not depends only on $L(\A)$ as follows:

\begin{cor}
Assume that $\A$ is free. Then $\rho^H$ is surjective if and only if 
the $b_2$-equality holds for $(\A,H)$. Namely, if $\A$ is free, then the surjectivity of 
$\rho^H$ is combinatorial for all $H \in \A$. 
\label{freesurj2}
\end{cor}

\noindent
\textbf{Proof}. Immediate by Theorem \ref{surj}.\owari
\medskip

Let us see how 
Corollary \ref{delfree} works in the following example:

\begin{example}
Let $\A$ be the Edelman-Reiner's example as in Example \ref{ERex1}. 
which was shown to be free with exponents $(1,5,5,5,5)$ in \cite{ER}. 
This is the first example of a free arrangement that has non-free restrictions. In this case, for 
$H:x_1-x_2-x_3-x_4-x_5=0$, 
$\A^H$ is not free, and the $b_2$-equality holds by Example \ref{ERex1}. Also, by Theorem \ref{SPOG}, 
$\pd \A'=1$. 
Let s detemine, by using results in this 
article, the structure of $D(\A^H)$ by numerical computations.

First, Theorem \ref{surj} shows that both $\pi^H$ and $\rho^H$ are 
surjective, and Theorem \ref{pdrest}  (5) 
shows that $\pd \A^H=\pd \A'=1$. Since $\A$ is free and $\rho^H$ is surjective, $D(\A^H)$ is generated by 
the Euler derivations together with $5$-derivations of degree $5$ since $\exp(\A)
=(1,5,5,5,5)$. This is the minimal set of generators since $\A^H$ is not free, and there is the unique relation among them. Let us 
determine the degree of 
the unique relation which completely determines the free resolution of $D(\A^H)$.
Let $\theta_E,\theta_2,\ldots,\theta_5$ be the minimal set of generators for $D(\A^H)$ with 
$\deg \theta_i=5$ for $i=2,3,4,5$.  By Theorem \ref{b2gen}, we may assume that 
$\theta_2,\ldots,
\theta_5 \in D(\A^H,m^H)$, and $Q'\theta_E,\theta_2,\ldots,
\theta_5$ form a generator for $D(\A^H,m^H)$, where 
$$
Q':=Q(\A^H,m^H)/Q(\A^H) \in S^H_5.
$$
Since $\A$ is free, Theorem \ref{Ziegler2} show that 
$\theta_2,\ldots,\theta_5$ are free basis for $D(\A^H)$. thus 
$$
Q'\theta_E=\sum_{i=2}^5 f_i \theta_i
$$
for some $f_i$, so the relation exists 
in degree $6$. As a 
consequence, we have the free resolution 
$$
0 \rightarrow S^H[-6] \rightarrow S^H[-1] \oplus (\oplus S^H[-5])^{\oplus 4} \rightarrow D(\A^H) \rightarrow 0.
$$
By Theorem \ref{SPOG}
again, we know that 
$$
0 \rightarrow S[-6] \rightarrow S[-1] \oplus (\oplus S[-5])^{\oplus 4} \rightarrow D(\A') \rightarrow 0.
$$
\label{ERex}
\end{example}

For the surjectivity of $\pi^H$, we have the following implication on the projective 
dimensions:

\begin{theorem}
Assume that $\pi^H$ is surjective. Then 
$$\pd \A' \le \pd \A$$
unless $\pd \A=0$. If $\pd \A=0$, then 
$\pd \A' \le 1$.
\label{pipd}
\end{theorem}

\noindent
\textbf{Proof}. 
Let $\theta_E,\theta_2,
\ldots,\theta_s$ be a set of generators for $D(\A)$ 
such that $\theta_i \in D_H(\A)$ for 
$i=2,\ldots,s$. Let 
$$
\theta_E^H:=\rho^H(\frac{Q(\A')}{Q(\A^H)}) \rho^H(\theta_E) 
\in D(\A^H,m^H)_{d},
$$
where $d:=|\A'|-|\A^H|+1$. Since $\pi^H$ is surjective, 
$$
\theta_E^H=\sum_{i=2}^s \pi^H(f_i \theta_i)=
\sum_{i=2}^s \rho^H(f_i \theta_i)
$$
for some $f_2,\ldots,f_s \in S$. 
Thus Proposition \ref{ER} shows that there is $\varphi \in D(\A')_{d-1}$ such that 
$$
\alpha_H \varphi=
\frac{Q(\A')}{Q(\A^H)}\theta_E-\sum_{i=2}^s f_i \theta_i.
$$
Since $\theta_i(\alpha_H)=0$ and $
\alpha_H \nmid \frac{Q(\A')}{Q(\A^H)}$, it holds that 
$\alpha_H \nmid \varphi(\alpha_H)$, i.e., 
$\varphi \not \in D(\A)$. Since $|\A'|-|\A^H|=
d-1=\deg \varphi$, Theorem \ref{B} shows that, for any 
$\theta \in D(\A')$, there is $f \in S$ such that 
$\theta -f \varphi \in D(\A)$. Thus 
$D(\A')$ is generated by $D(\A)$ together with $\varphi$. 
Since the relation among $D(\A)$ and $\varphi$ is 
$$
\alpha_H \varphi=
\frac{Q(\A')}{Q(\A^H)}\theta_E-\sum_{i=2}^s f_i \theta_i.
$$
which is of length one, $\pd \A'$ does not increase if 
$\pd \A \ge 2$. When $\A$ is free, 
see Theorem \ref{SPOG}. \owari
\medskip

\subsection{Application to Sylvester-Gallai theorem}
The surjectivity of $\rho^H$ has applications to the existence of double points. 
We say that $X \in L_2(\A)$ is \textbf{complete intersection} if 
$\A_X$ consists of two hyperplanes. Whem $\ell=3$ this is nothing but the double points, and 
starting from the well-known Sylvester-Gallai theorem and the Dirac-Motzkin conjecture, several 
researches exist. For example, see \cite{A8}. In \cite{A8}, it is proved that 
$\A$ has $X \in L_2(\A)$ of complete intersection on $H$ if $\A \setminus \{H\}$ is free when $\ell=3$ as follows:

\begin{prop}[\cite{A8}, Theorem 3.2]
Let $\A$ be an arrangement in $\K^3$, where $\K$ is a field of characteristic zero.
Let 
$mdr(\A):=\min\{n\mid D_H(\A)_n \neq (0)\}$. 
Then for $H \in \A$ with $|\A^H|>mdr(\A)$, there is 
$X \in L_2(\A)$ of complete intersection on $H$.
\label{dp}
\end{prop}

From the viewpoint of the surjectivity, we can generalize this as follows:

\begin{theorem}
Assume that $\rho^H$ 
is surjective. Then there is $X \in L_2(\A)$ of complete intersection on $H$. 
\label{surjdouble}
\end{theorem}

\noindent
\textbf{Proof}. 
By Lemma 
\ref{surjloc}, 
it suffices to show that there is a 
complete 
intersection codimension two flat of $\A_X$
on $H$ for some $X \in L_2(\A^H)$. So we may assume that $\A$ is an arrangement in $\K^3$. 
Let $r:=mdr(\A)$. Since $\exp(\A^H)=(1,|\A^H|-1)$, 
the surjectivity of $\rho^H$ implies that 
$r \le |\A^H|-1$. Thus Proposition \ref{dp} completes the proof. \owari
\medskip

Thus Proposition \ref{dp} can be regarded, by applying Theorem \ref{freesurj}, 
as a corollary of Theorem \ref{surjdouble}. Also, 
we can 
give a generalization of Theorem 1.7 
in \cite{A8}.

\begin{theorem}
Assume that $\A'$ is free, or the $b_2$-equality holds for $(\A,H)$. Then 
there is a complete intersection 
condimension two flat in $L(\A)$ on $H$.
\label{freedouble}
\end{theorem}

\noindent
\textbf{Proof}. Apply Theorems \ref{freesurj}, \ref{Y}, \ref{surjdouble} and Proposition \ref{b2local}.\owari
\medskip

\section{Further problems}

From the results in this article it is natural to ask the following question:

\begin{problem}
Does the surjectivity of $\rho^H$ depend 
only on $L(\A)$?
\label{surjcomb}
\end{problem}

The answer is NO by the author and Michael DiPasquale as follows, in the computation of which 
we used Macaulay2 in \cite{GS}:

\begin{example}[\cite{D}]
Let us recall the Ziegler's pair of arrangements in $\CC^3$ from \cite{Z2}. Namely, 
$\A_1$ is defined by
$$
xyz(x+y+z)(2x+y+z)(2x+3y+z)(2x+3y+4z)(x+ 3z)(x+ 2y+ 3z)=0,
$$
and $\A_2$ by 
$$
xyz(x+y+z)(2x+y+z)(2x+3y+z)(2x+3y+4z)(3x+5z)(3x+4y+5z)=0.
$$
They have isomorphic lattices, so they have the same combinatorics. However, it is known that they 
have different free resolutions. Namely, $D(\A_1)$ is generated by one degree five derivation and 
three degree six derivations, but $D(\A_2)$ is generated by four degree six derivations.  By this example, Ziegler proved that the minimal set of generators for logarithmic derivation modules are not combinatorial, thus neither is the free 
resolution.

Now let $H:x=0$. They can be regarded as the same hyperplane in the lattice 
isomorphism $L(\A_1) \simeq L(\A_2)$. Since $|\A_i^H|=6$, it holds that 
$\exp(\A_i^H)=(1,5)$. Thus the above shows that 
$$
\rho^H_2:D(\A_2) \rightarrow D(\A_2^H)
$$
is not surjective since $D_H(\A_2)_{5}=0$. Let us prove that 
$$
\rho^H_1:D(\A_1) \rightarrow D(\A_1^H)
$$
is surjective. Since $\dim_\CC D_H(\A_1)_5=1$, say $\theta \in D_H(\A_1)_5$, it suffices to show that 
$\rho_1^H(\theta) \neq 0$. Equivalently, there are no $\varphi \in D(\A_1 \setminus \{H\})_4$ 
such that $\alpha_H \varphi=\theta$. By the computation, we can see that 
every $\phi \in D(\A_1 \setminus \{H\})_4$ belongs to $S \theta_E$. Thus $\rho_1^H(\theta)\neq 0$ and 
$\rho_1^H(\theta_E)$ form a basis for $D(\A_1^H)$. Hence $\rho_1^H$ is surjective.

\label{noncombex}

\end{example}

Still we have a lot of problems related to contents in this article as follows:

\begin{problem}
(1)\,\,
Can we formualte the addition-deletion results for projective dimensions when both $\A$ and $\A'$ are 
not NMPD?

\noindent
(2)\,\,
Can we determine the freeness of $\A'$ when $\A$ is not free? 
By the results above, we know that $\rho^H$ is surjective, but also know that the $b_2$-equality is not sufficient.
At least we need that $\pd \A=\pd \A^H+1$ if $\pd \A$ is not maximal. 


\noindent
(3)\,\,
Does $\pd \A$ depend only on $L(\A)$?

\end{problem}

\end{document}